\documentclass{amsart}
\usepackage{graphicx} % Required for inserting images
\usepackage{amssymb,latexsym,amsmath,extarrows}
\usepackage{enumitem}
\usepackage{tikz}
\usepackage{xcolor}
\usepackage{hyperref,url}
\usepackage{cite}
\usepackage[T1]{fontenc}
\hypersetup{hidelinks}

\numberwithin{equation}{section}

\newtheorem{theorem}{Theorem}[section]
\newtheorem{lemma}[theorem]{Lemma}

\newtheorem{proposition}[theorem]{Proposition}
\newtheorem{remark}[theorem]{Remark}

\newtheorem{example}[theorem]{Example}

\newtheorem{definition}[theorem]{Definition}
\newtheorem{corollary}[theorem]{Corollary}

\newtheorem{conjecture}[theorem]{Conjecture}

\newcommand{\R}{\mathbb R}
\newcommand{\B}{\mathbb{B}}
\newcommand{\M}{\mathbb{U}}
\newcommand{\Pcal}{\mathbb{P}}
\newcommand{\N}{\mathcal N}
\newcommand{\one}{\mathbf 1}
\newcommand{\supp}{\operatorname{supp}}
\newcommand{\Inc}{I}
\newcommand{\Fourier}{\mathcal F}

\title[Falconer estimate for circular and cinematic Furstenberg problem]{Falconer estimate for circular and cinematic Furstenberg problem}

\author{Xiang Xie}
\address{Academy of Mathematics and Systems Science, Chinese Academy of Sciences, Beijing, China}
\email{xiexiang@amss.ac.cn}
\date{26 October 2026}

\begin{document}

\begin{abstract}
  We prove the Falconer-type bound \(\dim_{\mathrm H}F\geq u+1\) for planar circular \((u,v)\)-Furstenberg sets whenever \(0<u\leq1\), \(0<v\leq3\), and \(2u+v\geq3\). The same bound holds for smooth cinematic curve families. This resolves the circular Furstenberg conjecture in \cite{GreenHarrisOuRenTammen2025Circular} in its full Falconer range \(2u+v\geq3\).
\end{abstract}

\maketitle

\section{Introduction}

Given $0<u\leq1$ and $0\leq v\leq3$, a set $F\subset\R^2$ is called a
\emph{circular $(u,v)$-Furstenberg set} if there is a set
$E\subset\R^2\times(0,\infty)$ with $\dim_{\mathrm H}E\geq v$ such that
\[
  \dim_{\mathrm H}\bigl(F\cap C(x,t)\bigr)\geq u
  \qquad ((x,t)\in E),
\]
where $C(x,t)$ is the circle with center $x$ and radius $t$. Our main result
establishes Falconer type bound for circular Furstenberg set in the range $2u+v\geq3$.

\begin{theorem}\label{main_thm_Hausdorff_form}
Let $0<u\leq1$ and $0<v\leq3$. If $F\subset\R^2$ is a circular
$(u,v)$-Furstenberg set and $2u+v\geq3$, then
\[
  \dim_{\mathrm H}F\geq u+1.
\]
\end{theorem}

The exponent $u+1$ is sharp in this range. Theorem~\ref{main_thm_Hausdorff_form} proves the Falconer type circular Furstenberg estimate in \cite{GreenHarrisOuRenTammen2025Circular} at the range $2u+v\geq3$.

We also establish the same Falconer estimate for $C^\infty$ cinematic curves.

\begin{definition}[Smooth cinematic graph family]\label{def:cinematic-family}
Let $Z\Subset Z^+\subset\R^3$ and $I\Subset I^+\subset\R$ be compact
parameter and graph charts inside open sets, and let
$g\in C^\infty(Z^+\times I^+)$. Write $z=(z_1,z_2,z_3)$ and
\[
  \Sigma_z=\{(t,g(z,t)):t\in I\},\qquad z\in Z.
\]
We call $\{\Sigma_z:z\in Z\}$ a \emph{smooth cinematic family} parameterized by $Z$  if
\begin{equation}\label{eq:cinematic-nondegeneracy}
    \left|\det\partial_{(z_1,z_2)}(g,g_t)\right|\geq c_0,
\end{equation}
\begin{equation}\label{eq:cinematic-curvature}
    \left|\det\partial_z(g,g_t,g_{tt})\right|\geq c_0
 \quad\text{on }Z^+\times I^+
\end{equation}
for some $c_0>0$. The charts may be subdivided finitely, and all
constants below may depend on fixed derivative bounds for $g$ and $c_0$.
\end{definition}

\begin{remark}\label{rem:cinematic-localization}
We use the graph formulation of cinematic curves in \cite[Definition~2.7]{GreenHarrisOuRenTammen2025Circular}. Indeed, our result can be generalized to smooth level-set family with $\Phi:Z^+\times \R^2\to \R$ satisfying Phong-Stein rotational curvature condition (\cite{BeltranHickmanSogge2021Survey}) and Sogge's cinematic curvature condition in \cite{Sogge1991Cinematic,KolasaWolff1999Kakeya},

\[
  \Sigma_z=\{y\in\R^2:\Phi(z,y)=0\},\qquad z\in Z.
\]
More precisely, we use the localization argument given in \cite[Appendix~A, especially (A.5)]{PramanikYangZahl2022Furstenberg} to transfer the level-set family to the graph family $\{(t,g(z,t))\}$. Therefore, Phong-Stein rotational curvature condition becomes \ref{eq:cinematic-nondegeneracy} and Sogge's cinematic curvature condition becomes \ref{eq:cinematic-curvature} (see \cite{GreenHarrisOuRenTammen2025Circular}).
\end{remark}

\begin{definition}[Cinematic Furstenberg set]\label{def:cinematic-furstenberg}
For $0<u\leq1$ and $0\leq v\leq3$, a set $F\subset\R^2$ is a
\emph{cinematic $(u,v)$-Furstenberg set} for a smooth cinematic family
$\{\Sigma_z:z\in Z\}$ if some $E\subset Z$ with
$\dim_{\mathrm H}E\geq v$ satisfies
\[
  \dim_{\mathrm H}(F\cap\Sigma_z)\geq u\qquad(z\in E).
\]
\end{definition}

\begin{theorem}[Cinematic Falconer estimate]
\label{thm:cinematic-hausdorff}
Let $0<u\leq1$, $0<v\leq3$, and $2u+v\geq3$. Every cinematic
$(u,v)$-Furstenberg set for a smooth cinematic family as in
Definition~\ref{def:cinematic-family} satisfies
\[
  \dim_{\mathrm H}F\geq u+1.
\]
\end{theorem}
\subsection{Linear Furstenberg estimate}

The circular problem is a variant of the linear Furstenberg problem. A set
$F\subset\R^2$ is a linear $(u,v)$-Furstenberg set if there is a family
$\mathcal L$ of affine lines with $\dim_{\mathrm H}\mathcal L\geq v$ such
that $\dim_{\mathrm H}(F\cap\ell)\geq u$ for every $\ell\in\mathcal L$.
The classical conjecture, concerning a line in every direction, was formulated
by Furstenberg in unpublished work and recorded by Wolff
\cite{Wolff1999Kakeya}. A long
series of subsequent improvements culminated in the work of Orponen and
Shmerkin \cite{OrponenShmerkin2023ABC} and Ren and Wang \cite{RenWang2023Furstenberg}.

\begin{theorem}[sharp linear Furstenberg estimate \cite{OrponenShmerkin2023ABC,RenWang2023Furstenberg}]\label{thm:linear-furstenberg}
Let $0<u\leq1$ and $0\leq v\leq2$. Every linear $(u,v)$-Furstenberg set
$F\subset\R^2$ satisfies
\begin{equation}\label{eq:linear-furstenberg}
  \dim_{\mathrm H}F
  \geq
  \min\left\{u+v,\frac{3u+v}{2},u+1\right\}.
\end{equation}
\end{theorem}

The Furstenberg estimate is related to exceptional sets estimate for orthogonal projections in the plane by point--line duality. By analogy, we call $u+v$ the Kaufman-type
branch, which governs when $v\leq u$, and $u+1$ the Falconer-type branch, which governs when $u+v\geq2$.

The Kaufman type estimate ($u+v$, when $v\leq u$) can be established by $L^2$ method (see \cite{Wolff1999Kakeya}). The Falconer type estimate ($u+1$, when $u+v\geq2$) was proved by Fu and Ren in \cite{FuRen2024Incidence}. Their proof exploits high-low method developed by Guth, Solomon and Wang \cite{GuthSolomonWang2019Incidence}. 

In the intermediate range $u<v<2-u$, the middle term governs. For $v=1$,
Wolff established the lower bound
$\dim_{\mathrm H}F\geq\max\{2u,u+\tfrac12\}$
\cite{Wolff1999Kakeya}. Orponen and Shmerkin later obtained a strict
improvement over $2u$ when $u<\min\{v,1\}$, using multiscale analysis and
Bourgain's projection theorem \cite{orponen2023hausdorff}. They
subsequently proved \eqref{eq:linear-furstenberg} for Ahlfors-David regular line families. Ren and Wang treated
semi-well-spaced configurations by a high-low argument and combined
this case with the regular one through multiscale decomposition, proving
Theorem~\ref{thm:linear-furstenberg} in full generality
\cite{RenWang2023Furstenberg}. Recently, Orponen, Pyörälä and Yi generalize  \eqref{eq:linear-furstenberg} from line to transversal family \cite{orponen2025furstenbergsettheoremtransversal}.

\subsection{Circular Furstenberg problem}

More recently, Green, Harris, Ou, Ren, and Tammen proposed the following
circular Furstenberg conjecture
\cite{GreenHarrisOuRenTammen2025Circular}.

\begin{conjecture}[Circular Furstenberg conjecture]\label{conj:circular}
Let $0<u\leq1$ and $0\leq v\leq3$. Every circular $(u,v)$-Furstenberg set
$F\subset\R^2$ satisfies
\[
  \dim_{\mathrm H}F
  \geq
  \min\left\{u+v,\frac{5u+v}{3},u+1\right\}.
\]
\end{conjecture}

The history begins with Wolff's work on circular Kakeya sets. His estimates
imply \cite{Wolff2000LocalSmoothing}
\[
  \dim_{\mathrm H}F\geq1+v
  \qquad (u=1,\ 0\leq v\leq1).
\]
F\"assler, Liu, and Orponen established the
sharp Kaufman bound $u+v$ for $v\leq u$ \cite{FasslerLiuOrponen2024Circular}. Combined with Liu's work \cite{Liu2023Circular}, they give
\[
\dim_{\mathrm H}F\geq
u+\max\{v/3,\min\{u,v\}\}
\]

Zahl extended the Kaufman estimate to smooth curve families forbidding tangencies
of any fixed order $k$ \cite{Zahl2026MaximalCurves}. (When $k=2$, it is known as cinematic family curves).

The conjectured Falconer-type branch is $\dim_{\mathrm H}F\geq u+1$ when
$2u+v\geq3$. Green, Harris, Ou, Ren, and Tammen proved
\cite{GreenHarrisOuRenTammen2025Circular}
\begin{align}
  \dim_{\mathrm H}F&\geq\min\{4u+v-3,u+1\},\label{eq_first_bound}\\
  \dim_{\mathrm H}F&\geq\min\{2u+\tfrac{2v}{3}-1,u+1\}.
\end{align}
Consequently, the $u+1$ bound holds when $3u+v\geq4$ or
$3u+2v\geq6$. They also obtained the bound \eqref{eq_first_bound} for cinematic Furstenberg set. Theorem~\ref{main_thm_Hausdorff_form}  and Theorem \ref{thm:cinematic-hausdorff} extends it to the full
conjectured range $2u+v\geq3$.

\subsection{Connection with restricted projections}

Point--line duality connects linear Furstenberg sets with exceptional sets
for planar projections. Similarly, curved Furstenberg problems is related to restricted projection to line in $\R^3$ (see duality \cite{KaenmakiOrponenVenieri2025Marstrand,
GreenHarrisOuRenTammen2025Circular}). More precisely, the analogous connection involves sine waves. 

For $z=(z_1,z_2,z_3)\in\R^3$, define restricted projection
\[
  \gamma(\theta)=\frac1{\sqrt2}(\cos\theta,\sin\theta,1),
  \qquad \rho_\theta(z)=\langle z,\gamma(\theta)\rangle,
\]
and define the sine wave (which is also a cinematic curve)
\[
  \Gamma_z=\left\{\left(\theta,
    \frac{z_1\cos\theta+z_2\sin\theta+z_3}{\sqrt2}\right):
    0\leq\theta<2\pi\right\}.
\]

\begin{definition}[Sine-wave Furstenberg set]
A set $F\subset\R^2$ is a sine-wave $(u,v)$-Furstenberg set if there is a
set $E\subset\R^3$ with $\dim_{\mathrm H}E\geq v$ such that
$\dim_{\mathrm H}(F\cap\Gamma_z)\geq u$ for every $z\in E$.
\end{definition}

We recall the Falconer type exceptional-set estimates for the projections $\rho_\theta$. For $A\subset\R^3$, write
\[
  \Theta_s(A)=\{\theta\in[0,2\pi):
                 \dim_{\mathrm H}\rho_\theta(A)<s\}.
\]

\begin{theorem}[Falconer estimate
  \cite{GanGuthMaldague2024Exceptional}]\label{thm:ggm-projection}
Let $A\subset\R^3$ be Borel, with $\dim_{\mathrm H}A=v$. For
$0\leq s\leq1$,
\[
  \dim_{\mathrm H}\Theta_s(A)
  \leq\max\left\{0,1+\frac{s-v}{2}\right\}.
\]
In particular, if $v>1$, then
$\dim_{\mathrm H}\Theta_1(A)\leq(3-v)/2$.
\end{theorem}

\begin{remark}
    Pramanik, Yang, and Zahl proved Kaufman estimate for $\rho_\theta$ \cite{PramanikYangZahl2022Furstenberg}, i.e. for
$0\leq s<\min\{\dim_{\mathrm H}A,1\}$,
\[
  \dim_{\mathrm H}\Theta_s(A)\leq s.
\] They obtain their estimate through incidence bounds
for curves, in the spirit of Wolff's circular maximal-function method
\cite{Wolff1997Circles}. This provides
another approach to estimates related to the circular Furstenberg problem.
\end{remark}

The Falconer estimate has the following equivalent sine-wave formulation, with a common
parameter set on every curve.

\begin{corollary}
  \label{cor:ggm-sine-wave}
Let $A\subset\R^3$ and $\Theta\subset[0,2\pi)$ be Borel, with
$\dim_{\mathrm H}A=v>1$ and $\dim_{\mathrm H}\Theta=u>0$. Define
\[
  F_{A,\Theta}=\{(\theta,\rho_\theta(z)):
                        \theta\in\Theta,\ z\in A\}.
\]
Writing $(F_{A,\Theta})_\theta=\{t:(\theta,t)\in F_{A,\Theta}\}$, we have
\[
  \dim_{\mathrm H}\{\theta\in\Theta:
     \dim_{\mathrm H}(F_{A,\Theta})_\theta<1\}
  \leq\frac{3-v}{2}.
\]
If $2u+v>3$, then $F_{A,\Theta}$ is a sine-wave $(u,v)$-Furstenberg set
and $\dim_{\mathrm H}F_{A,\Theta}\geq u+1$.
\end{corollary}

The proof of Gan, Guth and Maldague follows the idea of high-low method and employs fractal small caps decoupling for cone in $\R^3$. However, the common-parameter hypothesis for $\Theta$ distinguishes this corollary from the
unrestricted sine-wave Furstenberg problem. Moreover, their proof relies on $u$-spacing condition stemming from the common-parameter hypothesis and thus cannot be generalized to the general case.

\subsection{The high--low method}

\begin{figure}[htbp]
  \centering
  \includegraphics[width=.74\textwidth]{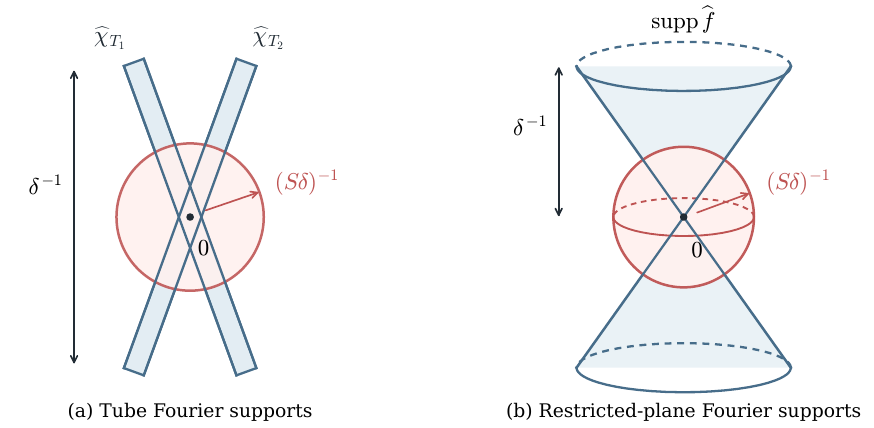}
  \caption{Fourier geometry of the high--low method. Left: Fourier
  rectangles associated with planar tubes. Right: Fourier tubes
  associated with planes of restricted normals lie near a cone.}
  \label{fig:high-low-frequency}
\end{figure}

The framework of our approach is the high--low method, introduced by Guth,
Solomon, and Wang \cite{GuthSolomonWang2019Incidence}. They study incidences
between $\delta$-tubes and $\delta$-balls in the plane. We briefly recall its idea in this part. Up to the natural
normalization and harmless enlargements, the incidence count is represented
by $\int fg$, where
\[
  f=\sum_{T\in\mathbb T}\chi_T,
  \qquad g=\sum_{B\in\mathbb B}\chi_B.
\]
Here $\chi_T$ and $\chi_B$ denote smooth approximations to characteristic
functions. The method analyzes this integral from Fourier perspective:
\[
  \int_{\R^2}fg
  =\int_{\R^2}\widehat f\,\overline{\widehat g}.
\]
For a tube $T$ of dimensions $\delta\times1$, the transform
$\widehat{\chi_T}$ is essentially supported in a rectangle of dimensions
$1\times\delta^{-1}$ centered at the origin, with its long axis normal to
$T$. Near $|\xi|\sim\delta^{-1}$, the Fourier rectangles corresponding to
$\delta$-separated directions have bounded overlap. More generally, outside
$B(0,(S\delta)^{-1})$, sufficient angular separation gives almost
orthogonality. This controls the high-frequency contribution through an
$L^2$ estimate; see Figure~\ref{fig:high-low-frequency}, left.

For the low-frequency contribution on $B(0,(S\delta)^{-1})$, let $\eta_{S\delta}$ be a smooth
approximate identity whose Fourier transform localizes to
$B(0,(S\delta)^{-1})$. The corresponding term takes the form
\[
  \int_{\R^2} f\,(g*\eta_{S\delta}).
\]
After harmless enlargements, this reduces to an incidence problem for
$S\delta$-balls and $S\delta$-tubes, which is handled by induction on
scales.

Gan, Guth, and Maldague \cite{GanGuthMaldague2024Exceptional} apply the same
principle to incidences between $\delta$-balls and localized
$\delta$-neighborhoods of planes with normals $\gamma(\theta)$. Again one
considers $\int fg$, with $f$ the sum of plane cutoffs and $g$ the sum of
ball cutoffs. A localized plane of dimensions $1\times1\times\delta$ has
Fourier transform essentially supported in a $1\times1\times\delta^{-1}$
tube through the origin, pointing in the direction $\gamma(\theta)$.
Consequently, $\widehat f$ is concentrated near the cone
\[
  \{\xi\in\R^3:\xi_1^2+\xi_2^2=\xi_3^2\}.
\]
The high-frequency part is controlled by fractal small-cap decoupling for
the truncated cone, while the low-frequency part is treated by induction;
see Figure~\ref{fig:high-low-frequency}.
For circular Furstenberg problem, Green et al. also introduced a high-low method for circles and balls 
\cite{GreenHarrisOuRenTammen2025Circular}. They express incidences through the
wave propagator. They use local smoothing, based on the sharp cone square
function estimate of Guth, Wang, and Zhang \cite{GuthWangZhang2020Cone},
to control high frequencies, and induction on scales to control low
frequencies. 

We conclude the philosophy of high-low method for Falconer type estimate:
\begin{enumerate}
    \item Fourier analysis estimate treats the high-frequency part
    \item Induction treats the low-frequency part.
\end{enumerate}

Compared with the previous work, we use refined decoupling for cone (Theorem \ref{thm:refined-decoupling}) to deal with the high-frequency part.

\subsection{Notation}

We write $A\lesssim B$ if $A\leq CB$ for an absolute constant $C$, and
$A\lessapprox B$ if, for every $\epsilon>0$,
$A\leq C_\epsilon\delta^{-\epsilon}B$. The relations $\gtrsim$, $\sim$,
$\gtrapprox$, and $\approx$ are defined analogously.
We write $A \lesssim_n B$ to emphasize the implicit constant may depend on $n$. When a scale is used in the article, we replace it
by a comparable dyadic scale without further comment.

The letters $C$ and $C_0$ denote fixed geometric constants, while $K$
denotes a regularity parameter, usually with $K\lessapprox1$.

\subsection{Acknowledgment}
The author is deeply grateful to Professor Chenjie Fan for his patient
guidance, generous support, and steady encouragement throughout the author's
mathematical studies. Professor Fan's mentorship has been invaluable in
shaping the author's path in mathematics.

\subsection{Funding Information and conflict of interest}
This work has been partial supported by NSFC Grant 12471232, 12288201, CAS Project for Young Scientists in Basic Research YSBR-031.

We declare that we have no conflict of interest in this work.

\subsection{AI Use Statement}
The circular-case result and its proof were developed by the author.
ChatGPT (OpenAI) assisted in developing the extension to smooth cinematic
families, including Theorem~\ref{thm:cinematic-hausdorff} and
Section~\ref{sec:cinematic-extension}, using the circular-case argument
as a starting point. It also assisted with language editing. The author
reviewed and revised the AI-assisted material and takes full
responsibility for the mathematical claims, proofs, citations, and
conclusions of this paper.

\section{Overview}

In this section, we outline the proof of main theorem \ref{main_thm_Hausdorff_form}.  For convenience, we simplify some definitions and give the informal version of key theorems; the precise
definitions and estimates appear in the subsequent sections. First, we formulate the relevant discretized circle-ball incidence problem.

\subsection{The discretized problem}

For a bounded set \(A\subset\R^d\), we use \(|A|_\delta\) in this overview
for the least number of balls of radius \(\delta\) needed to cover \(A\).
In the precise formulation below, \(|A|_\delta\) counts the dyadic
\(\delta\)-cubes meeting \(A\). The two quantities are comparable up to
a constant.

\begin{definition}[$(\delta, s, C)$-set]
     \textit{For $\delta \in 2^{-\mathbb{N}}$, $s \in [0, d]$ and $C > 0$, a non-empty bounded set $A \subset \mathbb{R}^d$ is called a $(\delta, s, C)$-set if}
\begin{equation*}
    |A \cap B(x, r)|_{\delta} \le C r^s |A|_{\delta}, \quad \forall x \in \mathbb{R}^d, r \in [\delta, 1].
\end{equation*}
 where we define $|A|_\delta$ as the dyadic $\delta$-covering number of $A$,
 \[
 |A|_\delta:= \#\{Q \in \mathcal{D}_{\delta}(\mathbb{R}^d) :   Q \cap A \neq \emptyset\}.
 \]
A family of dyadic $\delta$-cubes is a $(\delta,s,C)$-set if its union
has this property. For tangent rectangles with dimension $\delta\times\delta^{1/2}$, the corresponding condition
is imposed on their centers at the longitudinal scale $\delta^{1/2}$.
\end{definition}

Let $\Pcal$ be a family of $\delta$-neighborhoods of circles.  If
\[
  p=\N_\delta(C(x,t)),
  \qquad x\in\R^2,\quad 1\leq t\leq 2,
\]
define its dual ball by
\[
  \widehat p:=\N_\delta(x,t)\subset\R^3,
  \qquad
  \widehat{\Pcal}:=\{\widehat p:p\in\Pcal\}.
\]

\begin{definition}
We call $\Pcal$ a $(\delta,s,C)$-set if $\widehat{\Pcal}$ is a
$(\delta,s,C)$-set in $\R^3$.
\end{definition}

The discretized estimate needed for the main theorem is the following.

\begin{theorem}[Informal version of Theorem \ref{Thm_discretized}]\label{thm:discrete}
Let $\Pcal$ be a $(\delta,v,C)$-set of $\delta$-neighborhoods of circles
whose radii lie in $[1,2]$, and let $\B$ be a family of $\delta$-balls in
$\R^2$.  For every $p\in\Pcal$, let $\B(p)\subset\B$ consist of balls
contained in $p$, and suppose that $\B(p)$ is a $(\delta,u,C)$-set.  If
$2u+v\geq3$, then
\[
  |\B|\gtrapprox \delta^{-(u+1)}.
\]
The implicit constant may depend on $u,v$, and $C$.
\end{theorem}

The precise formulation in Definitions~\ref{definition__discrete_ball}
and~\ref{definition_discrete_circle} uses dyadic
cubes.  That formulation is better suited to the multiscale structure of the
argument.

\subsection{Two-scale decomposition}

The main new ingredient is to split the circle--ball incidence problem at the
intermediate scale $\delta^{1/2}$.  For each $p\in\Pcal$, cover the shading
$\B(p)$ by $\delta\times\delta^{1/2}$ rectangles tangent to $p$, and denote
the resulting family by $\M(p)$.  Set
\[
  \M:=\bigcup_{p\in\Pcal}\M(p).
\]
For simplicity, we may assume that every $U\in\M$ contains
$\sim N$ balls of $\B$, every $U\in\M$ is associated with $\sim M$ circles
of $\Pcal$, and every $p\in\Pcal$ contains $\sim W$ rectangles of $\M$. This part can be rigorously formalized through a series of pigeonholing.

\begin{definition}
Define
\begin{align*}
  \Inc(\Pcal,\M)
    &:=\sum_{p\in\Pcal}\sum_{\substack{U\in\M\\ U\subset p}}1,\\
  \Inc(\M,\B)
    &:=\sum_{B\in\B}\sum_{\substack{U\in\M\\ B\cap U\neq\varnothing}}1,\\
  \Inc(\Pcal,\B)
    &:=\sum_{p\in\Pcal}\sum_{\substack{B\in\B\\ B\cap p\neq\varnothing}}1.
\end{align*}
% The precise restricted incidence relation used later will be denoted by
% $\Inc_U(\Pcal,\B)$ and defined in Section~3.
\end{definition}

The uniformity assumptions give the bookkeeping identities
\begin{equation}\label{eq:two-scale-bookkeeping}
  \Inc_{\M}(\Pcal,\B)
  \sim M\Inc(\M,\B)
   \sim N\Inc(\Pcal,\M)
  \sim MN|\M|\sim NW|\Pcal|.
\end{equation}
Here $\Inc_{\M}$ counts incidence triples $(p,U,B)$ with
$U\subset p$ and $B\cap U\neq\varnothing$, as formalized in
Definition~\ref{define_incidence}. These identities concern restricted
incidences, rather than all incidences between $\Pcal$ and $\B$. 

\begin{figure}[ht]
  \centering
  \includegraphics[width=.47\textwidth]{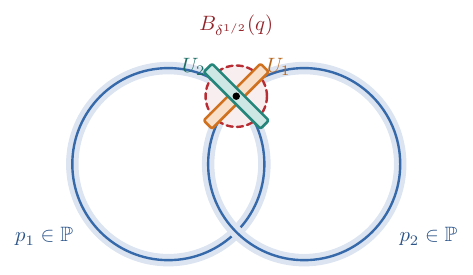}
  \caption{ $\delta^{1/2}$ scale decomposition of  circular incidence problem  }
  \label{fig:two-scale}
\end{figure}

\subsection{High--low estimates}
We will apply the high-low method to two incidence estimates $ \Inc(\M,\B)$ and $ \Inc(\Pcal,\M)$ respectively. 
We first see the incidence relation between
$\M$ and $\B$.

\begin{theorem}[Informal version of Theorem \ref{high_low_1}]\label{thm:hl-BM}
Let $S=\delta^{-\eta}$ with $0<\eta\ll1$.  At least one of the following
alternatives holds:
\begin{enumerate}[label=\textup{(\roman*)}]
  \item The high-frequency contribution dominates,
  \[
    \Inc(\M,\B)
      \lessapprox \delta^{-1/4}|\M|^{1/2}|\B|^{1/2}.
  \]
  Equivalently, under the uniformity assumptions above,
  \[
    N^2|\M|\lessapprox\delta^{-1/2}|\B|.
  \]
  \item The low-frequency contribution dominates,
  \[
    \Inc(\M,\B)\lessapprox \frac1S\Inc(\M,\B^S),
    \qquad
    \B^S:=\{\N_{S\delta}(B):B\in\B\}.
  \]
\end{enumerate}
\end{theorem}

After localization to $\delta^{1/2}$-balls and dilation by
$\delta^{-1/2}$, this becomes a standard incidence problem for
$\delta^{1/2}$-tubes and $\delta^{1/2}$-balls. The alternatives follow
from the planar high--low inequality in Proposition 2.1 
\cite{GuthSolomonWang2019Incidence}.

We next treat the incidence relation between $\Pcal$ and $\M$.  The relevant
geometry becomes more transparent after passing to dual objects in $\R^3$.

\begin{lemma}[Informal version of Lemma \ref{dual_incidence_relation}]\label{lem:duality}
There is a correspondence between $\delta$-circles and $\delta$-balls in
$\R^3$, and between $\delta\times\delta^{1/2}$ rectangles and
$\delta\times\delta^{1/2}\times1$ planks with restricted directions, such
that
\[
  U\subset p
  \quad\Longleftrightarrow\quad
  \widehat p\subset\widehat U.
\]
More explicitly, if the long axis of $U$ has direction
$(\cos\theta,\sin\theta)$, then the principal directions of $\widehat U$ are
\[
  (\sin\theta,-\cos\theta,1),\qquad
  (\cos\theta,\sin\theta,0),\qquad
  (-\sin\theta,\cos\theta,1),
\]
with corresponding side lengths $1$, $\delta^{1/2}$, and $\delta$; see figure \ref{fig:cone-duality}.
\end{lemma}

We are thus led to incidences between balls and planks with restricted
directions. The corresponding high--low alternatives are as follows.

\begin{theorem}[Informal version of \ref{thm_high_low_2}]\label{thm:hl-PM}
With $S=\delta^{-\eta}$, at least one of the following alternatives holds:
\begin{enumerate}[label=\textup{(\roman*)$'$}]
  \item The high-frequency contribution dominates,
  \[
    W^4|\Pcal|\lessapprox |\M|\delta^{-3/2}.
  \]
  \item The low-frequency contribution dominates,
  \[
    \Inc(\M,\Pcal)=\Inc(\widehat\M,\widehat\Pcal)
    \lessapprox \frac1S
      \Inc(\widehat\M,\widehat\Pcal^{\,S}),
  \]
  where
  $\widehat\Pcal^{\,S}:=\{\N_{S\delta}(\widehat p):
  \widehat p\in\widehat\Pcal\}$.
\end{enumerate}
\end{theorem}

A new ingredient of the proof is the application of refined decoupling theorem for cone. We give a sketch of proof.

\begin{proof}[Proof sketch]
Let
\[
  f:=\sum_{\widehat U\in\widehat\M}\one_{\widehat U},
  \qquad
  g:=\sum_{\widehat p\in\widehat\Pcal}\one_{\widehat p}.
\]
Partition $S^1$ into arcs $\omega$ of length $\delta^{1/2}$.  For
$\theta\in S^1$, set
\[
  \mathbf u_\theta=(\sin\theta,-\cos\theta,1),\quad
  \mathbf v_\theta=(\cos\theta,\sin\theta,0),\quad
  \mathbf w_\theta=(-\sin\theta,\cos\theta,1)
\]
and
\[
  \tau_\theta
  :=\left\{
    \xi_1\mathbf u_\theta+
    \xi_2\mathbf v_\theta+
    \xi_3\mathbf w_\theta:
    |\xi_1|\leq1,\ |\xi_2|\leq\delta^{-1/2},\
    |\xi_3|\leq\delta^{-1}
  \right\}.
\]
The boxes $\tau_\theta$ are comparable whenever
$|\theta_1-\theta_2|\leq\delta^{1/2}$, so we write $\tau_\omega$ for any
representative with $\theta\in\omega$. 

Let $\M_\omega$ be the rectangles whose long axes point in a direction
$(\cos\theta,\sin\theta)$ with $\theta\in\omega$, and let
$f_\omega$ be the corresponding part of $f$.  Then
$\widehat{f_\omega}$ is essentially supported in $\tau_\omega$.  Put
\[
  \tau_\omega^h:=\tau_\omega\setminus
    B(0,S^{-1}\delta^{-1}),
  \qquad
  f^h:=\Fourier^{-1}\!\left[
    \widehat f\,\one_{\{|\xi|\geq S^{-1}\delta^{-1}\}}
  \right],
\]
and define $f^\ell$ using the complementary low-frequency cutoff.  The union
of the high-frequency boxes $\tau_\omega^h$ forms a bounded-overlap cover of a
unit neighborhood of a truncated cone. See Figure \ref{fig:cone-duality}.

\begin{figure}[ht]
  \centering
  \includegraphics[width=.52\textwidth]{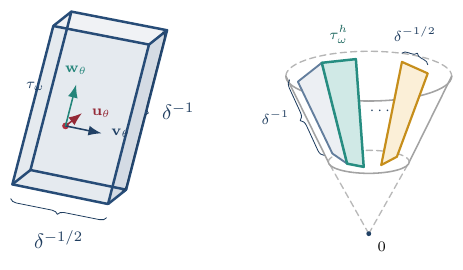}
  \caption{A frequency plank and high-frequency sectors on
  the truncated cone.}
  \label{fig:cone-duality}
\end{figure}

Suppose first that the high-frequency contribution dominates.  Up to harmless
constants,
\[
  \Inc(\widehat\Pcal,\widehat\M)
  \lesssim \delta^{-3}\int f^h g
  \leq \delta^{-3}\|f^h\|_{L^6(X)}\|g\|_{L^{6/5}},
\]
where $X:=\bigcup_{\widehat p\in\widehat\Pcal}\widehat p$ is a union of
$\delta$-balls, each meeting $\sim W$ planks from $\widehat\M$.  The refined
decoupling estimate for the cone gives
\[
  \|f^h\|_{L^6(X)}
  \lessapprox
  W^{1/3}
  \left(\sum_\omega\|f_\omega^h\|_6^6\right)^{1/6}
  \sim W^{1/3}(|\widehat\M|\delta^{3/2})^{1/6},
\]
while
\[
  \|g\|_{L^{6/5}}
  \sim (|\widehat\Pcal|\delta^3)^{5/6},
  \qquad
  \Inc(\widehat\Pcal,\widehat\M)\sim W|\Pcal|.
\]
Combining these relations yields
\[
  W^4|\Pcal|\lessapprox|\M|\delta^{-3/2}.
\]
If the low-frequency contribution dominates, enlarging the dual
$\delta$-balls to radius $S\delta$ gives alternative \textup{(ii)$'$}.
\end{proof}

If both high-frequency alternatives hold, then
\[
  W^4N^2|\Pcal|\lessapprox\delta^{-2}|\B|.
\]
Because each $\B(p)$ is a $(\delta,u,C)$-set, regularization and dyadic
pigeonholing allow us to assume
\[
  W\gtrsim|\M(p)|\gtrsim\delta^{-u/2},
  \qquad
  WN\gtrsim|\B(p)|\gtrsim\delta^{-u}.
\]
Consequently,
\[
  |\B|\gtrapprox \delta^{-3u+2}|\Pcal|
  \gtrsim\delta^{-3u-v+2}
  \gtrsim\delta^{-u-1},
\]
where the last inequality uses $2u+v>3$.

\subsection{Induction on scales}

If either low-frequency alternative holds, then the original incidence count
is controlled by its counterpart at the coarser scale $S\delta$; schematically,
\[
  \Inc(\Pcal,\B)
  \lessapprox \frac1S\Inc(\Pcal^S,\B^S).
\]
Following the multiscale analysis in \cite{orponen2023hausdorff,RenWang2023Furstenberg} , we show that
$\Pcal^S$ and $\B^S$ still form an essentially regular configuration at scale
$S\delta$: the circle family has dimension $v$, and every circle carries a
shading of dimension $u$.  This structure permits an induction on scales.

\subsection{Refined decoupling for the cone}

For completeness, we record the refined decoupling input used above.  Let
\[
  \Gamma:=\{(z\cos\theta,z\sin\theta,z):
  \theta\in[0,2\pi),\ z\in[1/2,1]\}
\]
be the truncated cone in $\R^3$.

\begin{theorem}[Refined decoupling for the cone \cite{GanGuoGuthEtAl2026RestrictedPlanes}]\label{thm:refined-decoupling}
Let $R\geq1$, and suppose that
$\supp\widehat f\subset\N_{R^{-1}}(\Gamma)$.  Decompose
$\N_{R^{-1}}(\Gamma)$ into $\sim R^{1/2}$ essentially disjoint
$R^{-1}\times R^{-1/2}\times1$ sectors $\tau$, and write the associated wave
packet decomposition as
\[
  f=\sum_\tau f_\tau
   =\sum_\tau\sum_{T\in\mathbb{T}_\tau}f_T.
\]
Here $\mathbb{T}_\tau$ is a family of planks dual to $\tau$, each of dimensions
$1\times R^{1/2}\times R$, and
$\mathbb{T}:=\bigcup_\tau\mathbb{T}_\tau$.  Assume that
\[
  \|f_T\|_2\sim\|f_{T'}\|_2
  \qquad (T,T'\in\mathbb{T}).
\]
For every $\varepsilon>0$, take $\eta>0$ sufficiently small in terms of
$\varepsilon$.
Let $Y$ be a disjoint union of unit balls in $B(0,R)$, each of which meets at
most $W$ dilated planks $R^\eta T$, where $T\in\mathbb{T}$.  Then, for
$2\leq p\leq6$,
\begin{equation}\label{eq:refined-decoupling}
  \|f\|_{L^p(Y)}
  \lesssim_{\varepsilon,\eta}
  R^\varepsilon
   W^{\frac12-\frac1p}
  \left(\sum_{T\in\mathbb{T}}\|f_T\|_{L^p}^p\right)^{1/p}.
\end{equation}
\end{theorem}

\section{Preliminary}

\subsection{Discretized problem}

We introduce the discretized circular Furstenberg problem and state the
discrete estimate underlying Theorem~\ref{main_thm_Hausdorff_form} and Theorem \ref{thm:cinematic-hausdorff}. The
definitions follow the linear framework of
\cite{RenWang2023Furstenberg}, with tangent rectangles added to record the
geometry at scale $\delta^{1/2}$.
%definition for distinct and comparable

\begin{definition}[$\delta$-discretized ball]
\label{definition__discrete_ball}
For $\delta=2^{-n}$, $n\in\mathbb Z$, let
\[
  \mathcal D_\delta(\R^d)
  :=\left\{\prod_{i=1}^d[a_i\delta,(a_i+1)\delta):
                 (a_1,\ldots,a_d)\in\mathbb Z^d\right\}.
\]
We call these dyadic cubes $\delta$-discretized balls. For dyadic
$0<\delta<\Delta\leq1$, every $B\in\mathcal D_\delta(\R^d)$ has a
unique ancestor $B_\Delta\in\mathcal D_\Delta(\R^d)$. If
$\mathbb B\subset\mathcal D_\delta(\R^d)$, set
\begin{equation}\label{eq:dyadic-ancestors}
  \mathcal D_\Delta(\mathbb B):=\{B_\Delta:B\in\mathbb B\},
  \qquad
  \mathbb B[B_\Delta]:=\{B\in\mathbb B:B\subset B_\Delta\}.
\end{equation}
We also write $\mathcal D^d_\delta$ for the $\delta$ cubes contained in $[0,1)^d$
\end{definition}

\begin{definition}[$\delta$-discretized circle]
\label{definition_discrete_circle}
For $\widehat p\in\mathcal D_\delta(\R^3)$ with positive radius
coordinate, define
\begin{equation}\label{eq:discretized-circle}
  p:=\bigcup_{(x,t)\in\widehat p}\N_\delta(C(x,t)).
\end{equation}
We regard $p$ as indexed by its parameter cube $\widehat p$ and write
$\widehat{\mathbb P}:=\{\widehat p:p\in\mathbb P\}$ for a family of
discretized circles. We call $\mathbb P$ a $(\delta,s,K)$-set, or
$(\delta,s,K)$-Frostman, if $\widehat{\mathbb P}$ is a
$(\delta,s,K)$-set in parameter space.

For dyadic $\delta<\Delta\leq1$, let $\Delta$ scale ancestor of $p$ be
\begin{equation}\label{eq:coarsened-circle}
  p^\Delta:=\bigcup_{(x,t)\in\widehat p_\Delta}
                       \N_\Delta(C(x,t)).
\end{equation}
We also write enlargement circle $Cp:=p^{C\delta}$. The notation
$\mathcal D_\Delta(\mathbb P)$ denotes the circles indexed by
$\mathcal D_\Delta(\widehat{\mathbb P})$, and
\[
  \mathbb P[p^\Delta]
  :=\{q\in\mathbb P:\widehat q\subset\widehat p_\Delta\}.
\]
For every $(x,t)\in\widehat p$, the Lipschitz dependence of a circle on
its center and radius gives
\begin{equation}\label{eq:circle-neighborhood-comparison}
  \N_\delta(C(x,t))\subset p\subset\N_{4\delta}(C(x,t)).
\end{equation}
\end{definition}

\begin{definition}[$\delta$-discretized cinematic curve]
\label{def:discretized-cinematic}
Fix the family and chart of Definition~\ref{def:cinematic-family}. For a
dyadic parameter cube $q\in\mathcal D_\delta(\R^3)$ contained
in the parameter chart $Z$, define its corresponding $\delta$- cinematic curve, 
\begin{equation}\label{eq:discretized-cinematic}
  p_q:=\bigcup_{z\in q}\N_\delta(\Sigma_z).
\end{equation}
A family $\mathbb P=\{p_q:q\in\widehat{\mathbb P}\}$ is a
$(\delta,s,K)$-set when its dual set
$\widehat{\mathbb P}$ form a $(\delta,s,K)$-set in parameter space.
For dyadic $\delta<\Delta\leq1$ whose ancestors remain in the chart,
set
\[
  p_q^\Delta:=\bigcup_{z\in q_\Delta}\N_\Delta(\Sigma_z),
  \qquad Cp_q:=\bigcup_{z\in q}\N_{C\delta}(\Sigma_z).
\]
The notation $\mathcal D_\Delta(\mathbb P)$ and $\mathbb P[p_q^\Delta]$
refers to the corresponding parameter ancestors, as for circles.
\end{definition}

\begin{lemma}[Cinematic geometry]\label{lem:cinematic-jet-geometry}
After a finite subdivision into sufficiently small parameter charts, the
maps
\[
  J_t:z\longmapsto\bigl(g(z,t),g_t(z,t),g_{tt}(z,t)\bigr)
\]
are uniformly bi-Lipschitz, for $t$ in a fixed graph interval. Hence
\begin{equation}\label{eq:cinematic-jet-metric}
  |z-z'|\sim_g\|g(z,\cdot)-g(z',\cdot)\|_{C^2(I)}.
\end{equation}
Therefore, there exists a geometric constant $C_g$, for $z\in q$,
\begin{equation}\label{eq:cinematic-neighborhood-comparison}
  \N_{C\delta}(\Sigma_z)\subset Cp_q\subset\N_{C_gC\delta}(\Sigma_z),
\end{equation}
\end{lemma}

Lemma \ref{lem:cinematic-jet-geometry} follows directly by the curvature condition \ref{eq:cinematic-curvature}, since $|D_zJ_t|\geq c_0$. 

    % shading definition generalized from Wang-Wu
    
\begin{definition}[Circular or cinematic Furstenberg shading]
Let $\mathbb P$ be a finite family of either $\delta$-discretized
circles with radii in $[1,2]$ or $\delta$-discretized cinematic curves in a fixed graph chart.
Let
$\mathbb B\subset\mathcal D_\delta(\R^2)$ be finite. A
$(\delta,C_0)$-shading model $(\Pcal,\B)_\delta$ consists of nonempty subfamilies
$\mathbb B(p)\subset\mathbb B$, $p\in\mathbb P$, such that
\[
  B\subset C_0p\quad(B\in\mathbb B(p)),
  \qquad
  \mathbb B=\bigcup_{p\in\mathbb P}\mathbb B(p).
\]
If $\mathbb P$ is a $(\delta,v,K)$-set and each $\mathbb B(p)$ is a
$(\delta,u,K)$-set, we call $(\mathbb B,\mathbb P)_\delta$ a
$(\delta,C_0,u,v,K)$ circular or cinematic Furstenberg shading,
according to the type of $\mathbb P$.
\end{definition}

The following is the discretized version of
Theorems~\ref{main_thm_Hausdorff_form} and
\ref{thm:cinematic-hausdorff}. We pass from this estimate to the
Hausdorff-dimensional statements in Section~\ref{section 7}.

\begin{theorem}[Discretized estimate]\label{Thm_discretized}
Fix either the circular family or one smooth cinematic graph family
as in Definition~\ref{def:cinematic-family}. Let $0<u\leq1$,
$0<v\leq3$, and $2u+v\geq3$. For every
$\epsilon>0$ and $C_0\geq1$, there exist $\eta>0$ and $\delta_0>0$
such that, for every dyadic $0<\delta<\delta_0$, every
$(\delta,C_0,u,v,\delta^{-\eta})$-circular or cinematic Furstenberg
shading $(\mathbb{P},\mathbb{B})_\delta$ satisfies
\begin{equation}\label{eq:discretized-falconer-bound}
  |\mathbb B|\gtrsim_{\epsilon,u,v,C_0}
                    \delta^{\epsilon-u-1}.
\end{equation}
For cinematic shadings, the constants also depend on the fixed graph function.
\end{theorem}

We also record the geometry of curves at scale $\delta^{1/2}$. On a graph chart
$\Sigma=\{(x,g(x))\}$ with uniformly bounded $C^2$ norm, Taylor's
theorem gives
\begin{equation}\label{eq:tangent-linearization}
  |g(x)-g(x_0)-g'(x_0)(x-x_0)|
  \lesssim|x-x_0|^2\lesssim\delta
  \qquad(|x-x_0|\lesssim\delta^{1/2}).
\end{equation}
Thus a $\delta$-neighborhood of such an arc is comparable, up to constant, to a tangent rectangle of dimensions
$\delta\times\delta^{1/2}$. First, we present the precise definitions for the terms "distinct" and "comparable".

\begin{definition}
For a rectangular box $V\subset\R^d$ with side lengths
$a_1,\ldots,a_d$, let $CV$ be its dilation by $C$ about its center.
Write $V^\Delta:=\N_\Delta(V)$. We call two rectangular boxes
$V_1,V_2$ $C_0$-comparable if
\[
  C_0^{-1}V_1\subset V_2\subset C_0V_1;
\]
otherwise they are $C_0$-distinct. Unless specified otherwise,
``distinct'' and ``comparable'' refer to a fixed absolute comparison
constant, which we take to be $10$.
\end{definition}

% \begin{remark}
% We define $p^{\Delta}$ and $B_\Delta$ for $\delta$ balls and circles previously. 
% \end{remark}

Let $\mathbb{U}$ be a set of distinct $\delta\times\delta^{1/2}$ rectangles. We adapt the definition of a configuration  from  \cite{orponen2023hausdorff}. Specifically, we introduce $\delta^{1/2}$ scale tangent boxes $\mathbb{U}$ within the shading model $(\mathbb{P},\mathbb{B})_\delta$.

\begin{definition}[Configuration]
A $(\delta,C_0)$-configuration $(\mathbb P,\mathbb U,\mathbb B)_\delta$
consists of a $(\delta,C_0)$-shading and nonempty subfamilies
$\mathbb U(p)\subset\mathbb U$ for each  $p\in \mathbb{P}$ such that
\[
  U\subset C_0p\quad(U\in\mathbb U(p)),
  \qquad
  \mathbb U=\bigcup_{p\in\mathbb P}\mathbb U(p),
\]
$C_0 \mathbb{U}(p)$ forms a cover of $\mathbb{B}(p)$. Moreover,
\begin{equation}\label{eq:balanced-tangent-shading}
  \#\{B\in\mathbb B(p):B\subset C_0U\}
  \sim_{C_0}\frac{|\mathbb B(p)|}{|\mathbb U(p)|}
  \qquad(U\in\mathbb U(p)).
\end{equation}
% The rectangles in $\mathbb U(p)$ are tangent to a representative circle
% for $p$ at the scale in \eqref{eq:tangent-linearization}. Their fixed
% enlargements have bounded overlap along that circle, with a bound
% depending only on the enlargement and $C_0$.

We call this a $(\delta,C_0,u,K,I,J)$-nice configuration if
$\mathbb B(p)$ is a $(\delta,u,K)$-set, $\mathbb U(p)$ form a $(\delta^{1/2},u,K)$-set, and
\[
  K^{-1}I\lesssim|\mathbb B(p)|\lesssim KI,
  \qquad
  K^{-1}J\lesssim|\mathbb U(p)|\lesssim KJ.
\]
\end{definition}

% add figures to illustrate, explain \sim_{K}

We next give definitions according to incidence between circles, points and tangent rectangles. Let $\mathbb{P}$ be a set of $\delta$ discretized circles or cinematic curves, $\mathbb{B}$ be a set of $\delta$ discretized balls in $\mathbb{R}^2$, and $\mathbb{U}$ be a set of distinct $\delta\times\delta^{1/2}$ rectangles.

\begin{definition}[Incidences]\label{define_incidence}
For $C\geq1$, define
\begin{align}
  I^C(\mathbb P,\mathbb B)
    &:=\sum_{p\in\mathbb P}
             \#\{B\in\mathbb B:B\subset Cp\},
       \label{eq:circle-ball-incidence}\\
  I^C(\mathbb P,\mathbb U)
    &:=\sum_{p\in\mathbb P}
             \#\{U\in\mathbb U:U\subset Cp\}.
       \label{eq:circle-box-incidence}
\end{align}
For a finite family $\mathbb T$ of convex bodies, set
\begin{equation}\label{eq:body-ball-incidence}
  I^C(\mathbb T,\mathbb B)
  :=\sum_{T\in\mathbb T}\#\{B\in\mathbb B:B\subset CT\}.
\end{equation}
The restricted incidence count of a shading is
\begin{equation}\label{eq:shading-incidence}
  \Pi(\mathbb P,\mathbb B):=\sum_{p\in\mathbb P}|\mathbb B(p)|.
\end{equation}
Finally, define the incidence count through $\mathbb U$ by
\begin{equation}\label{eq:incidence-through-boxes}
  I^C_{\mathbb U}(\mathbb P,\mathbb B)
  :=\sum_{p\in\mathbb P}
       \sum_{\substack{U\in\mathbb U\\U\subset Cp}}
           \#\{B\in\mathbb B:B\subset CU\}.
\end{equation}
\end{definition}

\begin{definition}[Refinement]
A $(\delta,C_0)$-configuration
$(\mathbb P',\mathbb U',\mathbb B')_\delta$ is a refinement of
$(\mathbb P,\mathbb U,\mathbb B)_\delta$ if
$\mathbb P'\subset\mathbb P$ with $ |\mathbb P'| \gtrapprox_\delta|\mathbb P|$;
$\mathbb B'(p)\subset\mathbb B(p)$, $\mathbb U'(p)\subset\mathbb U(p)$ for every $p\in\mathbb P'$, and
\begin{equation}\label{eq:refinement-cardinalities}
    |\mathbb B'(p)| \gtrapprox_\delta|\mathbb B(p)|,
    \qquad
    |\mathbb U'(p)|\gtrapprox_\delta|\mathbb U(p)|.
\end{equation}
A refinement of a $(\delta,C_0,u,K,I,J)$-nice configuration is
$(\delta,C_0,u,KK_0,I,J)$-nice for some
$K_0\lessapprox_\delta1$. This follows by restricting each normalized
Frostman bound and using \eqref{eq:refinement-cardinalities}.
\end{definition}

The next lemma turns preservation of the total shading incidence into a refinement.

\begin{lemma}[Refinement lemma]\label{refine_lemma}
Let $(\mathbb P,\mathbb U,\mathbb B)_\delta$ be a
$(\delta,C_0)$-configuration with
$|\mathbb B(p)|\sim I$ and $|\mathbb U(p)|\sim J$ for every
$p\in\mathbb P$. Suppose that
$(\mathbb B',\mathbb P')_\delta$ is a shading with
$\mathbb P'\subset\mathbb P$ and
$\mathbb B'(p)\subset\mathbb B(p)$ such that
\[
  \Pi(\mathbb P',\mathbb B')
  \gtrapprox_\delta\Pi(\mathbb P,\mathbb B).
\]
Then there is a refinement
$(\mathbb P'',\mathbb U'',\mathbb B'')_\delta$ of the original
configuration satisfying
$\mathbb P''\subset\mathbb P'$ and
$\mathbb B''(p)\subset\mathbb B'(p)$ for every $p\in\mathbb P''$.
\end{lemma}

\begin{proof}
Since $|\mathbb B'(p)|\lesssim I$, the incidence hypothesis allows us
to retain a subfamily $\mathbb P''\subset\mathbb P'$ with
\[
  |\mathbb P''|\gtrapprox_\delta|\mathbb P|,
  \qquad
  |\mathbb B'(p)|\gtrapprox_\delta I
  \quad(p\in\mathbb P'').
\]
For each retained $p$, dyadically pigeonhole the positive numbers
\[
  n_p(U):=\#\{B\in\mathbb B'(p):B\subset C_0U\},
  \qquad U\in\mathbb U(p).
\]
Since $C_0\M(p)$ form a cover of $\B(p)$, it gives a subfamily
$\mathbb U''(p)$ on which $n_p(U)\sim n_p$ and
\begin{equation}\label{eq:refinement-box-mass}
  n_p|\mathbb U''(p)|\gtrapprox_\delta|\mathbb B'(p)|.
\end{equation}
Let $\mathbb B''(p)$ consist of the balls of $\mathbb B'(p)$ contained
in at least one $C_0U$ with $U\in\mathbb U''(p)$. Bounded overlap and
\eqref{eq:refinement-box-mass} imply
$|\mathbb B''(p)|\gtrapprox_\delta I$, and the selected boxes satisfy
the uniform condition \eqref{eq:balanced-tangent-shading}. Moreover,
the uniform condition gives $n_p\lesssim_{C_0} I/J$, so
\[
  |\mathbb U''(p)|\gtrapprox_\delta J.
\]
Taking the unions over $p\in\mathbb P''$ produces the required
configuration.
\end{proof}
\subsection{Fourier Integral Operator (FIO)}

\begin{definition}[Phase of a cinematic function]
\label{defintion_phase}
Let \(g(z,t)\) be as in Definition~\ref{def:cinematic-family}, where
\(z=(x,\tau)\in\R^2\times\R\). On a compact annular wedge
\[
\mathcal A=\left\{\xi\in\R^2:\tfrac12\leq|\xi|\leq2,\ \xi_1/\xi_2\in I^+\right\},
\]
define the homogeneous phase function
\begin{equation}\label{eq:cinematic-phase}
\phi(z,\xi)=\xi_2g\!\left(z,\frac{\xi_1}{\xi_2}\right).
\end{equation}
\end{definition}

\begin{definition}[Conic FIO with H\"ormander conditions]
\label{def_FIO}
    Given a smooth, real-valued phase function $\phi(x,\tau;\xi) \in C^\infty(\R^{n+1}\times\R^n)$, homogeneous of degree one in
\(\xi\) and smooth away from $\xi = 0$, let $a \in S^0(\mathbb{R}^{n+1} \times \mathbb{R}^n)$ with compact support in $z=(x,\tau)\in\mathbb{R}^{n}\times \R$. For $\lambda\geq1$, the associated rescaled conic FIO is
\begin{equation}\label{eq:cinematic-fio}
T_\lambda h(y)
=\int_{\R^2}
 e^{2\pi i\phi_\lambda(y,\xi)}
 a(y/\lambda,\xi)h(\xi)\,d\xi,
\qquad
\phi_\lambda(y,\xi)=\lambda\phi(y/\lambda,\xi).
\end{equation}
    We call FIO $T_\lambda$ satisfy H\"onmander conditions if $\phi$ satisfies the following (H1) and (H2)
\begin{itemize}
    \item \textbf{(H1) Non-degeneracy condition.} 
    \[
    \operatorname{rank} \partial^2_{x\xi}\phi(x,\tau;\xi) = n \quad \forall (x,\tau;\xi) \in \operatorname{supp} a.
    \]
    \item \textbf{(H2) Positive definiteness condition.} 
    Consider the Gauss map $G : \operatorname{supp} a \to S^n$ by $G(z;\xi) := \frac{G_0(z;\xi)}{|G_0(z;\xi)|}$, where
    \[
    G_0(z;\xi) := \bigwedge_{j=1}^n \partial_{\xi_j} \partial_z \phi(z;\xi).
    \]
    Then for all $(z, \xi_0) \in \operatorname{supp} a$,
\begin{equation}\label{eq:gauss-curvature}
    \partial^2_{\xi\xi} \langle \partial_z \phi(z,\xi), G(z,\xi_0) \rangle |_{\xi=\xi_0}
\end{equation}
    has rank $n-1$ with $n-1$ positive eigenvalues.
\end{itemize}
\end{definition}

For the phase function $\phi$ in \ref{eq:cinematic-phase} given by a cinematic function $g$, direct calculation shows that nondegenerate condition  \ref{eq:cinematic-nondegeneracy} implies \textbf{(H1)} and curvature condition \ref{eq:cinematic-curvature} implies ~\eqref{eq:gauss-curvature} has rank one, slightly weaker than \textbf{(H2)}. Since $n=2$, its nonzero eigenvalue has a fixed sign on each connected chart. After subdividing the charts and, if necessary, reversing the \(\tau\)-coordinate, we may assume this sign is positive. Hence the phase can be locally reduced to one satisfying \textbf{(H2)}.

Fix \(0<\epsilon\ll1\), \(1\leq r\leq\lambda^{1-\epsilon}\), and set
\(\nu=\epsilon^{1000}\). We describe the scale \(r\) wave-packet
decomposition for \(T_\lambda\), following
\cite[Section~3.3]{GanHeLiWu2025LocalSmoothing}.

Choose a smooth partition of unity \(\sum_\theta\psi_\theta=1\) on
\(\mathcal A\), subordinate to angular caps
\begin{equation}
  \label{eq_fio_frequency_cap}
  \theta
=\left\{\xi\in\mathcal A:
\left|\frac{\xi_1}{\xi_2}-t_\theta\right|
\leq r^{-1/2}\right\},
\qquad
t_\theta\in r^{-1/2}\mathbb Z.
\end{equation}
Each \(\theta\) is comparable to a rectangle of dimensions
\(1\times r^{-1/2}\). For each \(\theta\), let \(\mathcal V_\theta\) be a bounded-overlap
covering of the physical space \(B_r^2\) by dual rectangles
\(V\) of dimensions \(1\times r^{1/2}\). Specifically, for $V\in\mathcal{V}_\theta$,
\begin{equation}
    \label{eq_FIO_physical_box}
    V=\{(S,Y):|S-S_V|\leq r^{1/2},|Y-Y_V+t_\theta S|\leq 1\}
    \quad S_V\in r^{1/2}\mathbb Z ,Y_V\in\mathbb Z
\end{equation}
Choose smooth functions
\(\eta_V\) adapted to these rectangles such that
\(\sum_{V\in\mathcal V_\theta}\eta_V=1\).

For \(T=(\theta,V)\), \(V\in\mathcal V_\theta\), define 
\[
h_T
:=\mathcal F\!\left[
  \eta_V\,\mathcal F^{-1}(\psi_\theta h)
\right]
=\widehat{\eta_V}*(\psi_\theta h).
\]
The corresponding wave packet function is \(T_\lambda h_T\). In a r radius ball $B^3_r$, we have
\begin{equation}
    T_\lambda h =\sum_T T_\lambda h_T + RapDec(r)\|h\|_{L^2}
\end{equation}
We also denote $\theta_T=\theta,V_T=V$.

\begin{definition}[curved planks]
Assume that $\phi(0,\xi)=0$. In $B_r$ center at the origin, the local inverse \(\gamma_\lambda(u,s;\xi)\) is determined by
\[
\partial_\xi\phi_\lambda
 \bigl((\gamma_\lambda(u,s;\xi),s),\xi\bigr)=u.
\]
By the non-degeneracy condition \textbf{(H1)}, after a finite localization of the chart and frequency wedge, the inverse is well-defined.

Let \(\xi_\theta\) be the center of \(\theta\), and let
\(\widetilde V\) denote the anisotropic enlargement of \(V\in\mathcal{V}_\theta\) with
dimensions \(r^{2\nu}\times r^{1/2+\nu}\). The enlarged curved plank
associated with \(h_T\) is
\[
T=T_{\theta,V}
=
\left\{
\bigl(\gamma_\lambda(u,s;\xi_\theta),s\bigr):
u\in\widetilde V,\ |s|\leq r
\right\}.
\]
For general \(\phi,B_r=B_r(y_0)\), the definitions are applied in translated
coordinates to
\[
\phi_{\lambda,y_0}(y,\xi)
=\phi_\lambda(y_0+y,\xi)-\phi_\lambda(y_0,\xi),
\qquad
h^{y_0}(\xi)=e^{2\pi i\phi_\lambda(y_0,\xi)}h(\xi).
\]
\end{definition}

Since \( \partial_\xi \phi (z,\xi)=(g_t(z,t),g-tg_t)\), where $t=\xi_1/\xi_2$. Combining \ref{eq_fio_frequency_cap} and \ref{eq_FIO_physical_box}, direct calculation implies that
\begin{equation}\label{eq_shape_of_wavepacket}
    T_{\theta,V}=\{y=\lambda z\in B_r:|g_t(z,t_\theta)-s_V|\leq r^{1/2+\nu}/\lambda,|g(z,t_\theta)-y_V|\leq r^{2\nu}/\lambda\}
\end{equation}
where $ s_V:=S_V/\lambda,y_v:=Y_V/\lambda$ . Here we assume that $\phi(0,\cdot)=0$. For general $\phi$ and $B(\lambda b,r)$, set $s_V^b:= s_V+g_t(b,t_\theta)$ and $ y_V^b:=y_V+g(b,t_\theta)$, then
\[
 T_{\theta,V}=\{y=\lambda z\in B(\lambda b,r):|g_t(z,t_\theta)-s^b_V|\leq r^{1/2+\nu}/\lambda,|g(z,t_\theta)-y^b_V|\leq r^{2\nu}/\lambda\}
\]

\begin{proposition}[Wave-packet decomposition]
\label{prop_wave_packet}
Let \(\mathbb T\) be the wave packets in \(B^3_r\), constructed
above. We have
\begin{enumerate}
    \item $L^2$ orthogonality
\[
\operatorname{supp}h_T\subset C\theta_T,\qquad
\sum_{T\in\mathbb T}\|h_T\|_2^2\lesssim\|h\|_2^2,
\]
\item Let $T$ be the associated curved plank. Then $h_T$ is essentially supported on $T$. More precisely, for every \(N>0\),
\[\|T_\lambda h_T\|_{L^\infty(B_r\setminus T)}
\lesssim_N r^{-N}\|h\|_2.
\]

\item  By the fixed-\(s\)
\(L^2\) estimate and local constancy property, we have\[
\sup_{|s|\leq r}
 \|T_\lambda h_T(\cdot,s)\|_{L_x^2}
 \lesssim\|h_T\|_2,
\qquad
\|T_\lambda h_T\|_{L^\infty(B_r)}
 \lesssim r^{-1/4}\|h_T\|_2,
\]
Furthermore, for a rapidly decreasing weight \(w_{B_r}\) adapted to \(B_r\),
\begin{equation}\label{eq:cinematic-single-packet}
\|T_\lambda h_T\|_{L^6(w_{B_r})}\lesssim\|h_T\|_2.
\end{equation}
\end{enumerate}
\end{proposition}

We refer to \cite[Section~3.3]{GanHeLiWu2025LocalSmoothing} and \cite[Theorem~2.1]{BeltranHickmanSogge2021Survey} for detailed proofs of these properties.

The Fourier integral operators can be viewed as generalizations of  Fourier extension operator on the cone in Euclidean space. We state the standard cone phase, which is also corresponding to the circular case.

\begin{example}[The cone and the shape of a packet]
For the cone phase
\[
\phi_{\mathrm{cone}}(x,\tau;\xi)=x\cdot\xi+\tau|\xi|,
\qquad
\gamma_{\mathrm{cone}}(u,\tau;\xi)
=u-\tau\frac{\xi}{|\xi|},
\]
the set
\[
T_{\theta,V}^{\mathrm{cone}}
=
\left\{
\left(u-s\frac{\xi_\theta}{|\xi_\theta|},s\right):
u\in\widetilde V,\ |s|\leq r
\right\}
\]
is a straight plank. Before enlargement its dimensions are
\(r\times r^{1/2}\times1\). For a general phase,
\(s\mapsto\gamma_\lambda(u,s;\xi_\theta)\) may curve, while each
\(s\)-slice remains the image of the same dual boxes.
\end{example}

Next, we introduce the variable coefficient version of Theorem \ref{thm:refined-decoupling} \cite[Appendix~A, Theorem~A.3]{GanHeLiWu2025LocalSmoothing}.

\begin{theorem}[Refined decoupling for FIOs]
\label{thm:cinematic-refined-decoupling}
Fix \(0<\epsilon\ll1\), and
\(1\leq r\leq\lambda^{1-\epsilon}\). Let \(\mathbb T\) be a collection of
the enlarged packets in $B_r$ above such that
\[
\|T_\lambda h_T\|_{L^6(w_{B_r})}\sim A
\qquad(T\in\mathbb T).
\]
Put \(K_{\mathrm{wp}}=r^{\epsilon^{50}}\). Suppose \(Y\subset B_r\)
is a union of disjoint balls of radius \(K_{\mathrm{wp}}\), each
meeting at most \(W\) packets \(T\in\mathbb T\). Then
\begin{equation}\label{eq:cinematic-refined-decoupling}
\|T_\lambda h\|_{L^6(Y)}
\lesssim_{\epsilon}
r^\epsilon
\left(\frac{W}{|\mathbb T|}\right)^{1/3}
\left(
\sum_{T\in\mathbb T}
\|T_\lambda h_T\|_{L^6(w_{B_r})}^{2}
\right)^{1/2}
+\operatorname{RapDec}(r)\|h\|_2.
\end{equation}
This is the \(n=2\), \(p=6\) case of
\cite[Appendix~A, Theorem~A.3]{GanHeLiWu2025LocalSmoothing}.
Equivalently, by comparability of the output norms,
\begin{equation}\label{eq:cinematic-refined-sixth-sum}
\|T_\lambda h\|_{L^6(Y)}
\lesssim_{\epsilon}
r^\epsilon W^{1/3}
\left(
\sum_{T\in\mathbb T}
\|T_\lambda h_T\|_{L^6(w_{B_r})}^{6}
\right)^{1/6}
+\operatorname{RapDec}(r)\|h\|_2.
\end{equation}
\end{theorem}

\section{Multiscale structures of configurations}

For dyadic $\delta<\Delta\leq1$, we construct a configuration at scale
$\Delta$ from one at scale $\delta$. This provides the structure needed
for the low-frequency alternative. The construction follows the
multiscale approach of \cite{orponen2023hausdorff} and its formulation in
\cite{RenWang2023Furstenberg}. We adapt linear version to curves.

\begin{definition}
Let
\[
  \delta=\Delta_N<\Delta_{N-1}<\cdots<\Delta_1<\Delta_0=1
\]
be dyadic scales. A finite family
$\mathcal A\subset\mathcal D_\delta(\R^d)$ is
$\{\Delta_j\}_{j=0}^N$-uniform if there are positive dyadic numbers
$N_j$ such that
\begin{equation}\label{eq:multiscale-uniformity}
  \#\{Q'\in\mathcal D_{\Delta_j}(\mathcal A):Q'\subset Q\}
  \in[N_j,2N_j)
\end{equation}
for every $j=1,\ldots,N$ and every
$Q\in\mathcal D_{\Delta_{j-1}}(\mathcal A)$. 

For a family of
discretized circles, uniformity refers to its parameter cubes.
\end{definition}

 Repeated
dyadic pigeonholing at a fixed number of scales retains a uniform
subfamily with only polylogarithmic loss (see Lemma 2.15 in \cite{OrponenShmerkin2023ABC}).

The next two propositions are the circle analogues of the coarse
configuration construction and lifting argument in
\cite[Proposition 3.6 and the proof of Proposition 4.2]{RenWang2023Furstenberg}.

\begin{proposition}\label{induction_scale_1}
Let $(\mathbb P_0,\mathbb U_0,\mathbb B_0)_\delta$ be a
$(\delta,C_0,u,K,I,J)$-nice configuration, and let
$\delta<\Delta<1$ be dyadic. Assume that \(\mathbb{P}_0\) is \(\{\delta, \Delta, 1\}\)-uniform. Then there is a
refinement $(\mathbb P_1,\mathbb U_1,\mathbb B_1)_\delta$ such that
\begin{enumerate}
\item
$|\mathcal D_\Delta(\mathbb P_1)|
 \approx_\delta|\mathcal D_\Delta(\mathbb P_0)|$ and
$|\mathbb P_1[p^\Delta]|
 \approx_\delta|\mathbb P_0[p^\Delta]|$ for every
$p^\Delta\in\mathcal D_\Delta(\mathbb P_1)$.
\item
Writing $\mathbb P_\Delta=\mathcal D_\Delta(\mathbb P_1)$, there is
an associated $(\Delta,C'_0,u,K_1K,i,j)$-nice configuration
$(\mathbb P_\Delta,\mathbb U_\Delta,\mathbb B_\Delta)_\Delta$,
where $K_1\lessapprox_\delta 1$ and $C'_0$ is a fixed constant depending on the cinematic function $g$ and $C_0$.
\end{enumerate}
\end{proposition}

\begin{proof}
    For each $p\in \mathbb{P}$, applying dyadic pigeonholing on finitely many arclength charts, there exist $\mathbb{B}_1(p) \subseteq \mathbb{B}(p)$ with $|\mathbb{B}_1(p)|\gtrapprox |\mathbb{B}(p)|$ such that
    
    $\mathbb{B}_1(p)$ is $\{\delta,\Delta,\delta^{1/2},\Delta^{1/2},1\}$ uniform. Let $\mathbb{U}_1(p)\subseteq\mathbb{U}(p)$ be the covering of $\mathbb{B}_1(p)$.  Assume that $ |\mathcal{D}_\Delta(\mathbb{B}_1(p))|\sim i_2(p)$, and each $B_\Delta \in \mathcal{D}_\Delta(\mathbb{B}_1(p))$
    \[
    \#\{Q\in\mathbb{B}_1(p):Q\subseteq B_\Delta\} \sim i_1(p)
    \]
    Pigeonholing again, there exist $\{\delta,\Delta,1\}$-uniform set $\mathbb{P}_1 \subseteq \mathbb{P}$ with  \(|\mathbb{P}_1| \gtrapprox |\mathbb{P}|\) and for each \(p \in \mathbb{P}_1\), \(i_2(p) \sim i_2,i_1(p)\sim i_1\). Explicitly, $i_1i_2\sim |\B_1(p)|\approx I$.

For each $p^\Delta$, pigeonhole the multiplicity of a coarse ball
among the shadings of $\mathbb P_1[p^\Delta]$. This gives a family
$\mathbb B_\Delta(p^\Delta)$ for any $B_\Delta$ in which
\begin{equation}\label{eq:coarse-ball-multiplicity}
\#\{p\in\mathbb P_1[p^\Delta]:
 B_\Delta\in\mathcal D_\Delta(\mathbb B_1(p))\}
\sim X(p^\Delta)
\end{equation}
and
$X(p^\Delta)|\mathbb B_\Delta(p^\Delta)|
 \approx_\delta |\mathbb{P}_1[p^\Delta]|i_2$. These balls lie in a fixed enlargement of $p^\Delta$.

Pigeonholing for $ \mathcal{D}_\Delta(\mathbb{P}_1)$, there exist $\mathbb{P}_\Delta\subseteq \mathcal{D}_\Delta(\mathbb{P}_1)$ such that $ |\mathbb{P}_\Delta|\gtrapprox |\mathcal{D}_\Delta(\mathbb{P}_1)|$ and each $p^\Delta \in \mathbb{P}_\Delta$, $X(p^\Delta)\sim X$ and $ |\mathbb{B}_\Delta(p^\Delta)|_{\Delta^{1/2}}\sim j$.
For notational simplicity, we also denote the ball set $ \bigcup_{p^\Delta \in \mathbb{P}_\Delta} \mathbb{P}_1[p^\Delta]$ as $\mathbb{P}_1$.

Following the proof of Proposition 4.1 in \cite{orponen2023hausdorff}, it follows that \(\mathbb{B}_\Delta(p^\Delta)\) is a \((\Delta, u, K_1K)\)-set for \(K_1 \lessapprox_\delta 1\). 

Let
\[
i  :=  \dfrac{|\mathbb{P}_1[p^\Delta]|\;i_2}{ X}
\]
And define $ \mathbb{U}_\Delta(p^\Delta)$ be a collection of distinct $\Delta\times\Delta^{1/2}$ rectangles covering $\mathbb{B}_\Delta(p^\Delta) $.  Consequently,  $|\mathbb{U}_\Delta(p^\Delta)|\sim j,|\mathbb{B}_\Delta(p^\Delta)|\sim i$. By \ref{eq:tangent-linearization}, $ U \subseteq C_0 p$ implies $ U_\Delta \subseteq C'_0 p^\Delta$ (In fact, we can choose $C_0' = 10C^2_0\|g\|_{C^2(I)}$). And by uniformity,  \(\mathbb{U}_\Delta(p^\Delta)\) is a \((\Delta^{1/2}, u, K_1K)\)-set. 

Let $\mathbb{U}_\Delta := \bigcup_{p^\Delta} \mathbb{U}_\Delta(p^\Delta)$ and $\mathbb{B}_\Delta:=  \bigcup_{p^\Delta} \mathbb{B}_\Delta(p^\Delta)$.

To ensure the rectangles in $\mathbb{U}_\Delta$ are pairwise distinct, we consider a maximal $10$- distinct subset $ \bar{\mathbb{U}}_\Delta \subseteq \mathbb{U}_\Delta $. Mapping each $U\in\mathbb{U}_\Delta$ to a comparable $\bar{U} \in \bar{\mathbb{U}}_\Delta$, we construct $\bar{\mathbb{U}}_\Delta(p)$. It follows that $(\mathbb{P}_\Delta,\bar{\mathbb{U}}_\Delta, \mathbb{B}_\Delta)_\Delta$ is a $(\delta,100C'_0)$ configuration with \(|\bar{\mathbb{U}}_\Delta(p^\Delta)| \sim |\mathbb{U}_\Delta(p^\Delta)| \sim j\). We absorb 100 into $C'_0$.

Relabeling $\bar{\mathbb{U}}_\Delta$ as $\mathbb{U}_\Delta$, we conclude that $(\mathbb{P}_\Delta,\mathbb{U}_\Delta,\mathbb{B}_\Delta)_\Delta$ is as desired.

% We can assume rectangles in $\mathbb{U}_\Delta $ are distinct with each other. Actually, let $ \bar{\mathbb{U}}_\Delta$ be the maximal $10$ distinct set in $\mathbb{U}_0$. And for $U\in\mathbb{U}_\Delta$  let $\bar{U}$ be the set in $\bar{\mathbb{U}}_\Delta $ comparable  with $U$. Define $ \bar{\mathbb{U}}_0(p):=\bigcup_{U\in\mathbb{U}_0(p)} \bar{U}$. It is easy to see $ (\mathbb{P}_\Delta,\bar{\mathbb{U}}_\Delta, \mathbb{B}_\Delta)_\Delta$ is a $(\delta,100C_0)$ configuration. And since $\mathbb{U}_\Delta(p^\Delta)$ are distinct, 
% \[
% |\bar{\mathbb{U}}_\Delta(p^\Delta)| \sim |\mathbb{U}_\Delta(p^\Delta)| \sim j
% \]
% For simplicity, we still write $\mathbb{U}_\Delta$ as $\bar{\mathbb{U}}_\Delta$.

% So $ (\mathbb{P}_\Delta,\mathbb{U}_\Delta,\mathbb{B}_\Delta)_\Delta$ is as desired.
\end{proof}
   
   \begin{proposition}
       \label{induction_scale_2}
      Continue with the setting in  Proposition \ref{induction_scale_1}. Let $(\mathbb{P}'_\Delta,\mathbb{U}'_\Delta,\mathbb{B}'_\Delta)_\Delta$ be a refinement of $(\mathbb{P}_\Delta,\mathbb{U}_\Delta,\mathbb{B}_\Delta)_\Delta$ . Then there is a
refinement $(\mathbb P_2,\mathbb U_2,\mathbb B_2)_\delta$ of
$(\mathbb P_1,\mathbb U_1,\mathbb B_1)_\delta$ such that  such that
      \begin{enumerate}
          \item $\mathcal{D}_\Delta[\mathbb{P}_2]\subseteq \mathbb{P}_\Delta'$ and $\mathcal{D}_\Delta[\mathbb{B}_2]\subseteq \mathbb{B}_\Delta'$
          \item For every $p\in\mathbb P_2$ and $U\in\mathbb U_2(p)$, there is a
$U_\Delta\in\mathbb U'_\Delta(p^\Delta)$ with
$U\subset C_*U_\Delta$, where $C_*$ depends only on the fixed
geometric constants.
      \end{enumerate}
 \end{proposition}
   \begin{proof}
      Recall \((\mathbb{P}_1, \mathbb{U}_1,\mathbb{B}_1)\) is the output of Proposition \ref{induction_scale_1} applied to \((\mathbb{P}_0, \mathbb{U}_0,\mathbb{B}_0)\). 
      Let \(\mathbb{P}' := \bigcup_{p^\Delta \in \mathbb{P}'_\Delta} \mathbb{P}_1[p^\Delta]\), and for each \(p \in \mathbb{P}_2\), let \(\mathbb{B}'(p) := \{B \in \mathbb{B}_1(p) : B_\Delta \in \mathbb{B}'_{\Delta}(p^\Delta)\}\). So $ (\mathbb{P}',\mathbb{B}')$ is a sub shading of $(\mathbb{P}_1,\mathbb{B}_1)$ such that $\mathcal{D}_\Delta[\mathbb{P}']\subseteq \mathbb{P}_\Delta'$ and $\mathcal{D}_\Delta[\mathbb{B}']\subseteq \mathbb{B}_\Delta'$

Since for each $p^\Delta$ and \(B_\Delta\in\mathbb{B}_\Delta(p^{\Delta})\),
\[
|\{p \in \mathbb{P}_1[p^\Delta]: B_\Delta \in \mathcal{D}_\Delta(\mathbb{B}_1(p))\}| \sim X
\]
and \(|\mathbb{B}_1(p)\cap B_\Delta| \sim i_1\), \( |\mathbb{B}'_\Delta(p^\Delta)|\approx |\mathbb{B}_\Delta(p^\Delta)|\approx i\). We calcualte $\Pi(\mathbb{P}',\mathbb{B}')$ as follows:
      \begin{align*}
          \Pi(\mathbb{P}',\mathbb{B}') & =  \sum_{p^\Delta \in \mathbb{P}'_{\Delta}} \sum_{B_\Delta\in \mathbb{B}'_\Delta(p^\Delta)} |\{(p, B) : p \in  \mathbb{P}_1[p^\Delta], B \in \mathbb{B}_1(p) \cap B_\Delta\}| \nonumber \\
          & \sim \sum_{p^\Delta \in \mathbb{P}'_{\Delta}} \sum_{B_\Delta\in \mathbb{B}'_\Delta(p^\Delta)} X\cdot i_1 \approx X\cdot i_1 i\cdot |\mathbb{P}'_\Delta| \\
      \end{align*}
      Recalling \( X \cdot i \approx_\Delta |\mathbb{P}_1[p^\Delta]|  i_2 ;\;i_1 i_2 \approx I\), thus $    \Pi(\mathbb{P}',\mathbb{B}') \approx |\mathbb{P}'|I\approx \Pi(\mathbb{P},\mathbb{B})$. 
      
      By Lemma \ref{refine_lemma}, there is a refinement  \((\mathbb{P}_2, \mathbb{U}_2,\mathbb{B}_2)_\delta \),  such that $ \mathbb{P}_2 \subseteq \mathbb{P}',  \mathbb{B}_2 \subseteq \mathbb{B}'$. It follows that $\mathcal{D}_\Delta[\mathbb{P}_2]\subseteq \mathbb{P}_\Delta'$ and $\mathcal{D}_\Delta[\mathbb{B}_2]\subseteq \mathbb{B}_\Delta'$.

Finally, since $ \mathbb{U}_2(p)$ is a $\delta^{1/2}$ covering of $\mathbb{B}_2(p) \subseteq \mathbb{B}'(p)$, then for each $U \in \mathbb{U}_2(p)$, there exists $U_\Delta \in \mathbb{U}'_\Delta(p^\Delta)$ such that $U \subseteq C_*U_\Delta$ for a fixed constant.

\end{proof}

\section{high-low methods for the circular case}

In this section, we will prove incidences estimates using high low method . We first recall the planar high--low estimate of
Guth, Solomon, and Wang \cite{GuthSolomonWang2019Incidence}, in the
form used by Ren and Wang \cite{RenWang2023Furstenberg}. We then
combine the high--low decomposition with refined cone decoupling to study
incidences between $\delta$-circles and
$\delta\times\delta^{1/2}$ tangent boxes.

\begin{theorem}[Proposition 2.6 in \cite{RenWang2023Furstenberg} ]\label{high_low_1}
    Let $0<\epsilon\ll1$, and set $S=\delta^{-\epsilon}$. Let
$\mathbb B$ be a family of distinct $\delta$-discretized balls and
$\mathbb T$ a family of distinct $\delta\times1$ tubes in a fixed
bounded subset of $\R^2$. Then at least one of the following holds:
\begin{enumerate}
\item The high-frequency dominate
\begin{equation}\label{eq:planar-high-incidence}
 I^C(\mathbb T,\mathbb B)
 \lesssim_{\epsilon,C}
 \bigl(S\delta^{-1}|\mathbb B||\mathbb T|\bigr)^{1/2}.
\end{equation}
\item The low-frequency estimate
\begin{equation}\label{eq:planar-low-incidence}
 I^C(\mathbb T,\mathbb B)
 \lesssim_{\epsilon,C}
 S^{-1+\epsilon}I_{\cap}^{C'}(\mathbb T,\mathbb B^{S\delta}),
\end{equation}
where $C'=C'(C)\geq1$ is fixed. Here
$\mathbb B^{S\delta}=\{B_{S\delta}:B\in\mathbb B\}$ is a multi-set counted
with multiplicity, and
\begin{equation}\label{eq:planar-enlarged-incidence}
 I_{\cap}^{C'}(\mathbb T,\mathbb B^{S\delta})
 :=\sum_{T\in\mathbb T}
   \#\{B_{S\delta}\in\mathbb B^{S\delta}:
                   B^{S\delta}\cap C'T\neq\varnothing\}.
\end{equation}
\end{enumerate}
\end{theorem}

We next focus on the incidence estimate between $\delta$ circles $p$ and tangent boxes $U$. First, we translate it to
ball--plank incidences in $\R^3$.

\begin{definition}[Duality of tangent boxes]
\label{definition_duality}
Let $U$ be a $\delta\times\delta^{1/2}$ rectangle with center
$(x_*,y_*)$ and long-axis direction $(\cos\theta,\sin\theta)$,
where $\theta\in[0,\pi)$. Put $\Theta(U)=\{\theta,\theta+\pi\}$,
and, for $\alpha\in\Theta(U)$, define
\begin{equation}\label{eq:dual-tangent-branch}
 \widehat U_\alpha
 :=\N_\delta\left\{
 (x_*+z\sin\phi,y_*-z\cos\phi,z):
 1\leq z\leq2,\ |\phi-\alpha|\leq\delta^{1/2}
 \right\}.
\end{equation}
Angular distances are interpreted modulo $2\pi$. The dual of $U$ is
\begin{equation}\label{eq:dual-tangent-box}
 \widehat U:=\widehat U_\theta\cup\widehat U_{\theta+\pi}.
\end{equation}

Each branch is comparable to a
$\delta\times\delta^{1/2}\times1$ plank. Its principal directions are
\[
 \mathbf u_\alpha=(\sin\alpha,-\cos\alpha,1),\qquad
 \mathbf v_\alpha=(\cos\alpha,\sin\alpha,0),\qquad
 \mathbf w_\alpha=(-\sin\alpha,\cos\alpha,1).
\]
% \begin{equation}\label{eq:dual-principal-directions}
%  \mathbf u_\alpha=(\sin\alpha,-\cos\alpha,1),\qquad
%  \mathbf v_\alpha=(\cos\alpha,\sin\alpha,0),\qquad
%  \mathbf w_\alpha=(-\sin\alpha,\cos\alpha,1).
% \end{equation}
with corresponding side lengths $1$, $\delta^{1/2}$, and $\delta$; see figure \ref{fig:cone-duality}.

We use these comparable rectangular planks when taking dilation of
$\widehat U$.
\end{definition}

\begin{lemma}[Dual incidence relation]\label{dual_incidence_relation}
There is an absolute constant $C_1$ such that, for every
$\delta$-discretized circle $p$, tangent box $U$, and $C\geq1$,
\begin{enumerate}
\item $U\subseteq Cp$ implies
      $\widehat p\subseteq C_1C\widehat U$;
\item $\widehat p\subseteq C\widehat U$ implies
      $U\subseteq C_1C^2p$.
\end{enumerate}
\end{lemma}
\begin{proof}
Choose $(a,b,r)\in\widehat p$, and write $q=(x_*,y_*)$ and
$v=(\cos\theta,\sin\theta)$. Since $1\leq r\leq2$, $x\in Cp$ implies
\begin{equation}\label{eq:circle-quadratic-neighborhood}
 |x-(a,b)|^2=r^2+O(C\delta).
\end{equation}
Suppose first that $U\subseteq Cp$. Applying
\eqref{eq:circle-quadratic-neighborhood} to the center and to the
two ends of the long axis of $U$ gives
\begin{equation}\label{eq:dual-tangency-test}
 |q-(a,b)|^2-r^2=O(C\delta),\qquad
 |\langle q-(a,b),v\rangle|\lesssim C\delta^{1/2}.
\end{equation}
Thus, we can write $(a,b,r)$ as
\[
 (a,b,r)
 =(x_*+r\sin\phi,y_*-r\cos\phi,r)+O(C\delta).
\]
for some $\alpha\in\Theta(U)$ and
$|\phi-\alpha|\lesssim C\delta^{1/2}$.

Equivalently, its coordinates relative to the corresponding dual
plank are $O(1)$, $O(C\delta^{1/2})$, and $O(C\delta)$.
The same bounds hold throughout $\widehat p$, proving the first
assertion.

Conversely, $(a,b,r)$ in $C\widehat U$ gives, for one of the two
branches,
\[
 |q-(a,b)|^2-r^2=O(C^2\delta),\qquad
 |\langle q-(a,b),v\rangle|\lesssim C\delta^{1/2}.
\]
Therefore, for $|t|\lesssim\delta^{1/2}$, 
$|q+tv-(a,b)|^2-r^2\lesssim O(C^2\delta)$.
The transverse width of $U$ contributes only $O(\delta)$, so
$U\subseteq C_1C^2p$.
\end{proof}

\begin{corollary}
Writing $\widehat{\mathbb U}=\{\widehat U:U\in\mathbb U\}$,
we have
\begin{align}
 I^C(\mathbb P,\mathbb U)
 &\leq I^{C_1C}(\widehat{\mathbb P},\widehat{\mathbb U}),
 \label{eq:dual-incidence-forward}\\
 I^C(\widehat{\mathbb P},\widehat{\mathbb U})
 &\leq I^{C_1C^2}(\mathbb P,\mathbb U).
 \label{eq:dual-incidence-backward}
\end{align}
\end{corollary}

\begin{theorem}\label{thm_high_low_2}
Let $\widehat{\mathbb P}$ be a family of distinct
$\delta$-discretized balls in $B^3(0,10)$, and let $\mathbb U$
be a family of distinct $\delta\times\delta^{1/2}$ rectangles
in $\R^2$. Let $0<\epsilon\ll1$, and set
$S=\delta^{-\epsilon}$. Suppose that every
$\widehat p\in\widehat{\mathbb P}$ is contained in
$\sim W$ members of $C\widehat{\mathbb U}$. Then at least
one of the following holds:
\begin{enumerate}
\item The high-frequency dominates
\begin{equation}\label{eq:dual-high-incidence}
 I^C(\widehat{\mathbb P},\widehat{\mathbb U})
 \lesssim_{\epsilon,C}
 S^{O(1)}W^{1/3}\delta^{-1/4}
 |\widehat{\mathbb P}|^{5/6}|\mathbb U|^{1/6}.
\end{equation}
Equivalently,
\begin{equation}\label{eq:dual-high-richness}
 W^4|\mathbb P|
 \lesssim_{\epsilon,C} S^{O(1)}\delta^{-3/2}|\mathbb U|.
\end{equation}
\item The low-frequency dominates
\begin{equation}\label{eq:dual-low-incidence}
 I^C(\widehat{\mathbb P},\widehat{\mathbb U})
 \lesssim_{\epsilon,C}
 S^{-1+\epsilon}
 I_{\cap}^{C'}(\widehat{\mathbb U},\widehat{\mathbb P}^{S\delta}),
\end{equation}
where $C'=C'(C)\geq1$ is fixed. Here
$\widehat{\mathbb P}^{S\delta}
 =\{\widehat p_{S\delta}:\widehat p\in\widehat{\mathbb P}\}$
is counted with multiplicity, and
\begin{equation}\label{eq:dual-enlarged-incidence}
 I_{\cap}^{C'}(\widehat{\mathbb U},\widehat{\mathbb P}^{S\delta})
 :=\sum_{U\in\mathbb U}
 \#\{\widehat p_{S\delta}\in\widehat{\mathbb P}^{S\delta}:
       \widehat p_{S\delta}\cap C'\widehat U\neq\varnothing\}.
\end{equation}
% This multiset differs from the dyadic coarse-scale family
% $\mathcal D_{S\delta}(\widehat{\mathbb P})$ in
% Definition~\ref{definition__discrete_ball}.
\end{enumerate}
\end{theorem}

\begin{proof}

Set \(\kappa=\epsilon/4\) and \(S_1=S^{1-\kappa}\). Let \(\mathcal V\) consist of two plank branches associated with each \(\widehat U\in \widehat\M\).

Partition \(S^1\) into arcs \(\omega\) of length \(\delta^{1/2}\). For \(\theta\in\omega\), put
\[
\mathbf u_\theta=(\sin\theta,-\cos\theta,1),\quad
\mathbf v_\theta=(\cos\theta,\sin\theta,0),\quad
\mathbf w_\theta=(-\sin\theta,\cos\theta,1),
\]
and let \(\tau_\omega\) be a box comparable to
\[
\left\{\xi_1\mathbf u_\theta+\xi_2\mathbf v_\theta+
\xi_3\mathbf w_\theta:
|\xi_1|\leq1,\ |\xi_2|\leq\delta^{-1/2},\
|\xi_3|\leq\delta^{-1}\right\}.
\]

Denote by \(\mathcal V_\omega\) the planks in \(\mathcal V\) associated with \(\omega\). Specifically, recalling Definition \ref{definition_duality},
\[
\mathcal{V}_\omega:= \{\hat{U}_\theta:U\in\mathbb{U},\ \theta\in\Theta(U)\cap \omega\}.
\]
For each \(V\in\mathcal V_\omega\), choose a nonnegative smooth function \(\chi_V\) adapted to \(V\), which is bounded below on \(CV\) and vanishes outside \(2CV\). It follows that \(\widehat{\chi_V}\) is essentially supported on \(\tau_\omega\). Specifically,
\begin{equation}
    \label{rapid_decay_V}
   \widehat{\chi_V}(\xi)\lesssim_{M,\eta} \delta^{M},\qquad \xi\notin \delta^{-\eta}\tau_{\omega}.
\end{equation}
Set \(R=\delta^{-1}\) and let \(\nu=\epsilon^{1000}\) be sufficiently small, \(A=R^\nu\). By the Schwartz property, it suffices to consider the Fourier area \(A\tau_\omega\).

Write
\[
f_\omega=\sum_{V\in\mathcal V_\omega}\chi_V,
\qquad
f:=\sum_\omega f_\omega.
\]
For \(\mathcal A\subseteq\widehat{\mathbb P}\), set
\[
g_{\mathcal A}=\sum_{\widehat p\in\mathcal A}\mathbf 1_{\widehat p}.
\]
Then
\[
I^C(\mathcal A,\widehat{\mathbb U})\sim W|\mathcal{A}|
\lesssim_C\delta^{-3}\int f g_{\mathcal A}.
\]

Choose a smooth bump function \(\phi_l\) which equals \(1\) on \(|\xi|\leq S_1^{-1}\delta^{-1}\) and is supported in \(|\xi|\leq 2S_1^{-1}\delta^{-1}\). Divide \(f\) into high and low frequency parts:
\[
f^\ell
=\Fourier^{-1}[\phi_l\widehat f],
\qquad f^h=f-f^\ell.
\]
The low-frequency kernel satisfies, for every \(M>0\),
\[
|K_\ell(x)|
\lesssim_M (S_1\delta)^{-3}
\left(1+\frac{|x|}{S_1\delta}\right)^{-M}.
\]
Similarly, define \(f_\omega^\ell,f_\omega^h, \chi^h_V\). Then the Fourier transform of \(f_\omega^h\) is essentially supported on \(\tau_\omega^h:=A\tau_\omega\setminus B(0,S_1^{-1}\delta^{-1})\), which is comparable to the plank
\[
\left\{\xi_1\mathbf u_\theta+\xi_2\mathbf v_\theta+
\xi_3\mathbf w_\theta:
|\xi_1|\leq A,\ |\xi_2|\leq A\delta^{-1/2},\ S_1^{-1}\delta^{-1}\leq
|\xi_3|\leq A\delta^{-1}\right\}.
\]

Fix a sufficiently large constant \(C'=C'(C)\), and set
\[
m_S(\widehat p)
=\#\{V\in\mathcal V:
\widehat p_{S\delta}\cap C'V\neq\varnothing\}.
\]
Then
\[
I_{\cap}^{C'}(\widehat{\mathbb U},\widehat{\mathbb P}^{S\delta})\sim \sum_{\widehat p\in\widehat{\mathbb P}}m_S(\widehat p).
\]

A plank meets an \(S\delta\)-ball in volume \(O(S^2\delta^3)\). Since \(S/S_1=S^\kappa\), the kernel tail outside the \(S\delta\) ball is negligible. Consequently, for every \(\mathcal A\subseteq\widehat{\mathbb P}\),
\begin{align*}
    \delta^{-3}\int f^\ell g_{\mathcal A} &= \delta^{-3}\int f (g_{\mathcal A}*K_l)\\
    &\lesssim_{\epsilon,C}
S^2S_1^{-3}\sum_{\widehat p\in\mathcal A}m_S(\widehat p)
+O(\delta^{100})\\
&=
S^{-1+3\kappa}\sum_{\widehat p\in\mathcal A}m_S(\widehat p)
+O(\delta^{100}).
\end{align*}

If \(I^C(\widehat{\mathbb P},\widehat{\mathbb U})
 \lesssim_{\epsilon,C}
 S^{-1+\epsilon}
 I_{\cap}^{C'}(\widehat{\mathbb U},\widehat{\mathbb P}^{S\delta})\), the low-frequency alternative follows. Otherwise, we obtain
 \[
 \sum_{\widehat p\in\widehat{\mathbb P}}m_S(\widehat p) \ll_{\epsilon,C}S^{1-\epsilon}W|\widehat{\mathbb P}|.
 \]
Discarding at most half of the balls, we obtain a set
\(\mathcal A\subseteq\widehat{\mathbb P}\) such that
\(|\mathcal A|\geq|\widehat{\mathbb P}|/2\) and
\begin{equation}
    \label{multiplicity_condition}
    m_S(\widehat p)\lesssim_{\epsilon,C}S^{1-\epsilon}W
\qquad(\widehat p\in\mathcal A).
\end{equation}
Moreover, the high frequency part dominates:
\begin{align}
    I^C(\mathcal A,\widehat{\mathbb U})&\sim W|\Pcal|
\lesssim_C\delta^{-3}\int f^h g_{\mathcal A}\\
&\leq \delta^{-3}
\|f^h\|_{L^6(X)}
\|g_{\mathcal A}\|_{6/5},
\end{align}
where \(X:=\bigcup_{\widehat p\in\mathcal{A}}\widehat p\).

The union of the high-frequency boxes \(\tau_\omega^h\) forms an \(O(S)\)-overlap cover of an \(O(S)\)-neighborhood of the following truncated cone; see also Figure \ref{fig:cone-duality}.
\[
 \Gamma:=\{(z\cos\theta,z\sin\theta,z):
  \theta\in[0,2\pi),\ z\in[S_1^{-1}\delta^{-1},A\delta^{-1}]\}.
\]
By the Schwartz property, we can decompose \(f^h=\sum_{\lambda\leq AR} f^\lambda+\mathrm{RapDec}(R)\) into \(O(\log R)\) dyadic radial pieces, and set
\[
 S_1^{-1}R\lesssim\lambda\lesssim A R,\qquad
 \Gamma_\lambda
 :=\{(z\cos\theta,z\sin\theta,z):
       \theta\in[0,2\pi),\ |z|\sim\lambda\}.
\]
Set \(\tau^\lambda_\omega:=\tau^h_\omega \cap \{z\sim\lambda\}\). Let \(f^\lambda_\omega,\chi_V^\lambda\) be the parts of \(f^h_\omega,\chi^h_V\) whose Fourier transforms are truncated to the shell \(\{|z|\sim \lambda\}\). Thus
\[
f^\lambda = \sum_\omega f_\omega^\lambda =\sum_{\omega}\sum_{V\in\mathcal{V}_\omega}\chi_V^\lambda.
\]
Fix \(\lambda\), and cover the parameter region by \(O(S^3)\) balls \(B=B(z_B,c/S)\). Put \(r=\lambda/S\). Under the change of variables \(y=\lambda(z-z_B)\), the Fourier support of \(f^\lambda\) becomes an \(O(r^{-1})\)-neighborhood of a unit truncated cone.

Perform the \(r\)-scale wave packet decomposition (in rescaled form) on each \(B\). The frequency space is partitioned into \(\Omega_r=\{\gamma_\theta:\theta \in r^{-1/2}\mathbb Z\cap [0,2\pi)\}\), where
\[
\gamma_\theta=\left\{\xi_1\mathbf u_\theta+\xi_2\mathbf v_\theta+
\xi_3\mathbf w_\theta:
|\xi_1|\leq S,\ |\xi_2|\leq \lambda r^{-1/2},
|\xi_3|\sim\lambda\right\}.
\]
Each wave packet \(T\in\mathbb T_{\gamma}\) is a translate of
\[
\left\{\xi_1\mathbf u_\theta+\xi_2\mathbf v_\theta+
\xi_3\mathbf w_\theta:
|\xi_1|\leq
S^{-1},\ |\xi_2|\leq \lambda^{-1} r^{1/2},
|\xi_3|\sim \lambda^{-1}\right\}.
\]

Denote by \(f^\lambda_T\) the wave packet function. For each \(T\in \mathbb T_\gamma\), since \(\delta^{1/2}\leq r^{-1/2}\), there are \(O(S)\) caps \(\omega\) such that \(\omega\subseteq2\gamma\). Since \(\chi_V\) is supported on \(2CV\), the function \(f^\lambda_T\) is negligible unless \(R^\nu T\) intersects the physical supports \(\{C'V\in C' \mathcal{V}_\omega:\omega\subseteq 2\gamma\}\). Denote by \(\mathbb T\) the set of all such wave packets.

Recall that \(V\in\mathcal{V}_\omega\) is comparable to a translate of
\[
\left\{\xi_1\mathbf u_\theta+\xi_2\mathbf v_\theta+
\xi_3\mathbf w_\theta:
|\xi_1|\leq1,\ |\xi_2|\leq\delta^{1/2},\
|\xi_3|\leq\delta\right\}
\]
for any \(\theta\in\omega\). Compared with $R^\nu T$, since
\[
 A\left(\frac{\sqrt r}{\lambda}
             +\frac{r^{-1/2}}S\right)=o(\sqrt\delta),
 \qquad
 A\left(\frac1\lambda+\sqrt\delta\,r^{-1/2}
             +S_1\delta\right)=o(S\delta),
 \quad A/S=o(1),
\]
It implies that \(R^\nu T\subseteq \mathcal{N}_{S\delta}(C'V)\) for some \(C'\) and \(V \in \mathcal{V}_\omega\) for $T\in \mathbb T$.

By \eqref{multiplicity_condition}, each \(\hat{p}\in\mathcal{A}\) intersects at most \(O(S)W\) sets \(\mathcal{N}_{S\delta}(CV)\), and thus at most \(S^{O(1)}W\) wave packets \(T\in\mathbb T\), since each \(\mathcal{N}_{S\delta}(CV)\) contains at most \(S^{O(1)}\) distinct \(T\). We also have
\[
|\mathbb T|\leq S^{O(1)}|\mathcal V|\lesssim S^{O(1)} |\mathbb U|.
\]
Applying Theorem~\ref{thm:refined-decoupling} at scale \(r\) yields
\[
 \|f^\lambda\|_{L^6(X\cap B)}
 \lesssim_\epsilon S^{O(1)}W^{1/3}
       \bigl(\sum_{T\in\mathbb T}\|f^\lambda_T\|_{L^6}^6\bigr)^{1/6}.
\]
It is straightforward that each \(T\) is contained in at most \(S^{O(1)}\) sets \(\mathcal{N}_{S\delta}(CV)\). Thus \(\|f^\lambda_T\|_{L^6} \leq S^{O(1)} \delta^{1/4}\), and
\[
\|f^\lambda\|_{L^6(X\cap B)}\lesssim_\epsilon S^{O(1)}W^{1/3}   \bigl(|\mathbb U|\delta^{3/2}\bigr)^{1/6}.
\]
Summing over \(\lambda\) and \(B\), we conclude that
\[
 \|f^h\|_{L^6(X)}
 \lesssim_\epsilon S^{O(1)}W^{1/3}
       \bigl(|\mathbb U|\delta^{3/2}\bigr)^{1/6}.
\]
Since \(\|g_{\mathcal A}\|_{6/5} \sim (|\Pcal|\delta^3)^{5/6}\),
\[
I^C(\hat{\Pcal},\hat{\mathbb{U}}) \lesssim_{\epsilon,C} S^{O(1)}W^{\frac{1}{3}}\delta^{-1/4}|\hat{\Pcal}|^{\frac{5}{6}}|\hat{\mathbb{U}}|^{\frac{1}{6}},
\]
equivalently
\[
W^4|\mathbb P|
\lesssim_{\epsilon,C}
S^{O(1)}\delta^{-3/2}|\mathbb U|,
\]
which is the high-frequency alternative.
\end{proof}

\section{ Proof of the main theorem for the circular case}

In this section, we prove Theorem \ref{Thm_discretized} for the circular case.
\begin{theorem}\label{Thm_configuration}
Let $0<\epsilon\ll1$, $C_0,C>10$, $0<u\leq1$, and $0<v\leq3$, with
$2u+v\geq3$. There is $\delta_0>0$ such that the following holds for every
dyadic $\delta\in(0,\delta_0)$ and $K\geq1$. Suppose
$(\mathbb P_0,\mathbb U_0,\mathbb B_0)_\delta$ is a circular
$(\delta,C_0,u,K,I,J)$-nice configuration and $\mathbb P_0$ is a
$(\delta,v,K)$-set. Then there are a refinement
$(\mathbb P',\mathbb U',\mathbb B')_\delta$ and a family
$\mathbb B_*\subseteq\mathcal D_\delta(\R^2)$ such that
$\mathbb B'\subseteq\mathbb B_*\subseteq\mathbb B_0$ and
\begin{equation}\label{eq:configuration-incidence-bound}
  I^C_{\mathbb U'}(\mathbb P',\mathbb B_*)
  \lesssim_{\epsilon,C_0,C}
  K^4\delta^{1-\epsilon}|\mathbb P'||\mathbb B_*|.
\end{equation}
\end{theorem}
We will prove Theorem \ref{Thm_configuration} by induction on $\delta$. We first explain how this estimate implies Theorem~\ref{Thm_discretized}.

\begin{proof}[Deduction of Theorem~\ref{Thm_discretized}]

Fix $0<\epsilon<1$ and set $\eta=\epsilon^2/100$. Let
$(\mathbb B,\mathbb P)_\delta$ be a
$(\delta,C_0,u,v,\delta^{-\eta})$-circular Furstenberg shading. We construct a nice configuration from it. 

For each $p\in\mathbb P$, choose a bounded-overlap cover of its shading
by distinct tangent rectangles of dimensions $\delta\times\delta^{1/2}$.
Write $\mathbb U(p)$ for this cover. For dyadic $L\geq1$, set
\begin{align*}
  \mathbb U^L(p)
  &:=\left\{U\in\mathbb U(p):
    \#\{B\in\mathbb B(p):B\subseteq C_0U\}\sim L\right\},\\
  \mathbb B^L(p)
  &:=\bigl\{B\in\mathbb B(p):
    B\subseteq C_0U\text{ for some }U\in\mathbb U^L(p)\bigr\}.
\end{align*}
Dyadic pigeonholing first in each shading and then in the circle family
gives common values $J,L,I:=JL$ and a subfamily $\mathbb P_0\subseteq\mathbb P$
such that
\[
  |\mathbb P_0|\gtrapprox|\mathbb P|,
  \qquad
  |\mathbb U^L(p)|\sim J,
  \qquad
  |\mathbb B^L(p)|\sim I\gtrapprox|\mathbb B(p)|
  \quad(p\in\mathbb P_0).
\]
Set $\mathbb B_0(p)=\mathbb B^L(p)$ and
$\mathbb U_0(p)=\mathbb U^L(p)$. 

Let $\mathbb{U}_0 := \bigcup_{p} \mathbb{U}_0(p)  $ and $\mathbb{B}_0 := \bigcup_{p} \mathbb{B}_0(p)  $. To ensure the rectangles in $\mathbb{U}_0$ are pairwise distinct, we consider a maximal $10$- distinct subset $ \bar{\mathbb{U}}_0 \subseteq \mathbb{U}_0 $. Mapping each $U\in\mathbb{U}_0$ to a comparable $\bar{U} \in \bar{\mathbb{U}}_0$, we construct $\bar{\mathbb{U}}_0(p)$. It is straightforward to verify that $(\mathbb{P}_0,\bar{\mathbb{U}}_0, \mathbb{B}_0)_\delta$ constitutes a $(\delta,100C_0)$-configuration.
Since $\mathbb{U}_0(p)$ are distinct, 
\[
|\bar{\mathbb{U}}_0(p)| \sim |\mathbb{U}_0(p)| \sim J
\]
Furthermore, for any $\bar{U}\in \bar{\mathbb{U}}_0(p)$, the dilated rectangle $100C_0 \bar{U}$ contains at most $C_0^{10}$ distinct rectangles from $\mathbb{U}_0(p)$. Consequently,
\[
L\lesssim\#\{B\in \mathbb{B}_0(p),B\subset 100C_0 \bar{U}\} \lesssim C_0^{10} L 
\]
For notational simplicity, we still write $\mathbb{U}_0$ as $\bar{\mathbb{U}}_0$, and absorb the factor $100$ into $C_0$.

Since $\mathbb{B}(p)$ is a $(\delta,u,\delta^{-\eta})$-set and $|\mathbb{B}_0(p)|\gtrapprox |\mathbb{B}(p)|$, it follow that $(\mathbb{P}_0,\mathbb{U}_0,\mathbb{B}_0)_\delta$ form a $(\delta,C_0,u,K_0 \delta^{-\eta},I,J)$ nice configuration with $ \mathbb{P}_0$ a $(\delta,v,K_0 \delta^{-\eta})$-set, where $K_0\lessapprox 1$ .

Apply Theorem \ref{Thm_configuration} to $(\mathbb{P}_0,\mathbb{U}_0,\mathbb{B}_0)_\delta$ with $C=C(C_0)$ large enough to contain the shading
cover, then 
 \[
    I^C_{\mathbb{U}'}(\mathbb{P}',\mathbb{B}_*) \lesssim _{\epsilon,C_0,C}  K^4 \delta^{1-\epsilon} |\mathbb{P}'||\mathbb{B}_*|
    \]
    Since
    \begin{align*}
    I^{C}_{\mathbb{U}'}(\mathbb{P}',\mathbb{B}_*) \geq I^{C_0}_{\mathbb{U}'}(\mathbb{P}',\mathbb{B}')\geq \sum_{p\in\mathbb{P}'} |\mathbb{B}'(p)|\gtrapprox \delta^{\eta-u}|\mathbb{P}'|
\end{align*}
Therefore, $|\mathbb{B}_0|\geq |\mathbb{B}_*|\gtrsim_{\epsilon,C_0} \delta^{\epsilon}\delta^{-u-1}$.
\end{proof}

Now we prove Theorem \ref{Thm_configuration}.

Fix $S:=\delta^{-\eta}$, where $\eta = \epsilon^2$. Let $ \Delta :=S\delta$.

\noindent\textbf{Step 1: Uniformization at scale $\Delta$.}

Let $\mathbb{B}_S:=\mathcal{D}_{S\delta}(  \mathbb{B}_0)$. For dyadic number $k>0$, define
$$ \mathbb{B}_{S,k}:=\{ B_S \in \mathbb{B}_S, |\mathbb{B}_0[B_S]|\sim k\}, \mathbb{B}_k := \bigcup_{B_S \in \mathbb{B}_{S,k}} \mathbb{B}_0[B_S]$$

Pigeonholing gives $k$ such that $ \Pi(\mathbb{P}_0,\mathbb{B}_k ) \gtrapprox \Pi(\mathbb{P}_0,\mathbb{B}_0)$. By  Lemma \ref{refine_lemma}, there exists a refinement $ (\mathbb{P}_1,\mathbb{U}_1,\mathbb{B}_1)_\delta$ with $ \mathbb{B}_1 \subseteq \mathbb{B}_k.$ Also, we can assume $\mathbb{P}_1$ is $\{\delta,\Delta,1\}$ -uniform.
  
%define -uniform.

\noindent\textbf{Step 2: Construction of $\Delta$-scale configuration.}

Applying Lemma \ref{induction_scale_1} to $(\mathbb{P}_1,\mathbb{U}_1,\mathbb{B}_1)_\delta$, we obtain a refinement $(\mathbb{P}'_1,\mathbb{U}'_1,\mathbb{B}'_1)_\delta$ and a $(\Delta,100C_0,u,KK_0,i,j)$-configuration $(\mathbb{P}_\Delta,\mathbb{U}_\Delta,\mathbb{B}_\Delta)_\Delta$ with $K_0\lessapprox 1$.

By induction on $S\delta$, there exists a refinement $(\mathbb{P}'_\Delta,\mathbb{U}'_\Delta,\mathbb{B}'_\Delta)_\Delta$ of $(\mathbb{P}_\Delta,\mathbb{U}_\Delta,\mathbb{B}_\Delta)_\Delta$, along with a set $\mathbb{B}_{\Delta,*}$ satisfying $\mathbb{B}'_\Delta \subseteq \mathbb{B}_{\Delta,*} \subseteq \mathbb{B}_\Delta$, such that
\begin{equation}\label{eq:coarse-inductive-estimate}
  I^{C}_{\mathbb{U}'_\Delta}(\mathbb{P}'_{\Delta},\mathbb{B}_{\Delta,*}) \lesssim_{\epsilon,C_0,C} K_0^4K^4 (\delta S)^{1-\epsilon} |\mathbb{P}'_\Delta||\mathbb{B}_{\Delta,*}|.
\end{equation}

Next, by Lemma \ref{induction_scale_2}, there exists a refinement $(\mathbb{P}_2,\mathbb{U}_2,\mathbb{B}_2)_\delta$ of $(\mathbb{P}'_1,\mathbb{U}'_1,\mathbb{B}'_1)_\delta$ such that
\begin{equation}\label{eq:lifted-coarse-supports}
  \mathcal D_\Delta(\mathbb P_2)\subseteq\mathbb P'_\Delta,
  \qquad
  \mathcal D_\Delta(\mathbb B_2)\subseteq\mathbb B'_\Delta.
\end{equation}
Furthermore, for each $U \in \mathbb{U}_2$, there exists $U_\Delta \in \mathbb{U}'_\Delta$ such that $U \subseteq C_* U_\Delta$ for some absolute constant $C_*$.

Finally, let $\mathbb{B}_*:=\bigcup_{B_S\in\mathbb{B}_{\Delta,*}} \mathbb{B}_0[B_S]$. Since $\mathbb{B}_{\Delta,*}\subseteq \mathbb{B}_\Delta\subseteq\mathcal{D}_\Delta(\mathbb{B}_1)$, the uniformity established in Step ~1 implies that $\mathcal{D}_\Delta(\mathbb{B}_*)=\mathbb{B}_{\Delta,*}$ and that $|\mathbb{B}_*[B_S]|$ is uniform up to a constant factor for each $B_S \in \mathbb{B}_{\Delta,*}$.

\noindent\textbf{Step 3: Uniformization of the two-scale incidences.}

For $U\in\mathbb U_2$, let
\[
  n(U):=\#\{B\in\mathbb B_*:B\subseteq CU\},
  \qquad
  m(U):=\#\{p\in\mathbb P_2:U\subseteq Cp\}.
\]
Pigeonholing these two quantities gives $M,N\geq1$ and
$\mathbb P_3\subseteq\mathbb P_2$ such that
$|\mathbb P_3|\gtrapprox|\mathbb P_2|$ and, for every $p\in\mathbb P_3$,
\[
  \mathbb U_3(p):=
  \{U\in\mathbb U_2(p):m(U)\sim M,\ n(U)\sim N\}
  \quad\text{satisfies}\quad
  |\mathbb U_3(p)|\gtrapprox|\mathbb U_2(p)|.
\]
Keep the balls of $\mathbb B_2(p)$ covered by these boxes, and denote
them by $\mathbb B_3(p)$.  Put
$\mathbb U_3=\bigcup_{p\in\mathbb P_3}\mathbb U_3(p)$. Thus $(\mathbb{P}_3,\mathbb{U}_3,\mathbb{B}_3)_\delta$ is a refinement.

Let $D_0:=C_1C^2$ be the fixed dilation from
Lemma~\ref{dual_incidence_relation}. 
A final pigeonholing gives $\mathbb P_4\subseteq\mathbb P_3$ and $W\geq1$
such that $|\mathbb P_4|\gtrapprox|\mathbb P_3|$ and
\begin{equation}\label{eq:uniform-dual-degree}
  \#\{U\in\mathbb U_3:\widehat p\subseteq D_0\widehat U\}
  \sim W\qquad(p\in\mathbb P_4).
\end{equation}
Then \( \#\{U\in\mathbb U_3: U\subseteq Cp\}
  \lesssim W\). In particular,
\begin{equation}\label{eq:restricted-incidence-upper-uniform}
  I^C_{\mathbb U_4}(\mathbb P_4,\mathbb B_*)
  \lesssim NW|\mathbb P_4|.
\end{equation}

Set
\[
  \mathbb U_4:=\bigcup_{p\in\mathbb P_4}\mathbb U_3(p),
  \qquad
  \mathbb B_4:=\bigcup_{p\in\mathbb P_4}\mathbb B_3(p).
\]
The output refinement will be
$(\mathbb P',\mathbb U',\mathbb B')_\delta
=(\mathbb P_4,\mathbb U_4,\mathbb B_4)_\delta$.

\noindent\textbf{Step 4: High-low estimates for Two scale structures}

Consider the incidence estimates between $\delta$ discretized balls $ \mathbb{B}_*$ and $\delta\times\delta^{1/2}$ rectangles $\mathbb{U}_3$, applying Theorem \ref{high_low_1} at scale $\delta^{1/2}$, one of the following holds, 
\begin{enumerate}[label=\textup{(\roman*)}]
  \item \label{high_1} The high-frequency estimate gives
  \[
    \Inc^C(\M_3,\B_*)
      \lesssim_{\epsilon,C} \delta^{-1/4-\epsilon^2}|\M_3|^{1/2}|\B_*|^{1/2}.
  \]
  Equivalently, under the uniformity assumptions above,
   \begin{equation}\label{eq:planar-high-uniform}
   N^2|\M_3|\lesssim_{\epsilon,C}\delta^{-1/2-2\epsilon^2}|\B_*|.
  \end{equation}
  \item \label{low_1} The low-frequency estimate gives
  \begin{equation}\label{eq:planar-low-uniform}
  I^C(\mathbb{U}_3,\mathbb{B}_*) \lesssim_{\epsilon,C} S^{-1+\epsilon^2} I_{\cap}^{C'}(\mathbb U_3,\mathbb B_*^{S\delta}).
  \end{equation}

\end{enumerate}

Next, apply Theorem~\ref{thm_high_low_2} to $\hat{\Pcal}_4$ and $\hat{\mathbb{U}}_3$ using
\eqref{eq:uniform-dual-degree}. At least one of the following holds:
\begin{enumerate}[label=\textup{(\roman*)$'$}]
  \item\label{high_2} The high-frequency estimate gives
  \begin{equation}\label{eq:dual-high-uniform}
    W^4|\mathbb P_4|
    \lesssim_{\epsilon,C}
   S^{O(1)} \delta^{-3/2}|\mathbb U_3|.
  \end{equation}
  \item\label{low_2} The low-frequency estimate gives
  \begin{equation}\label{eq:dual-low-uniform}
    I^{D_0}(\widehat{\mathbb P}_4,\widehat{\mathbb U}_3)
    \lesssim_{\epsilon,C}
    S^{-1+\epsilon^2}
    I_{\cap}^{C'}(\widehat{\mathbb U}_3,
                  \widehat{\mathbb P}_4^{S\delta}).
  \end{equation}
\end{enumerate}

\medskip
\noindent\textbf{Both high-frequency estimates hold.}
Combining \eqref{eq:planar-high-uniform} and
\eqref{eq:dual-high-uniform}, gives
\begin{equation}\label{eq:combined-high-estimates}
  W^4N^2|\mathbb P_4|
  \lesssim_{\epsilon,C}
  \delta^{O(\epsilon^2)}\delta^{-2}|\mathbb B_*|.
\end{equation}
The regularity of the refined shadings, their tangent covers, and the
parameter family gives
\begin{equation}\label{eq:regularity-degree-lower-bounds}
  W\gtrapprox K^{-1}\delta^{-u/2},
  \qquad
  WN\gtrapprox K^{-1}\delta^{-u},
  \qquad
  |\mathbb P_4|\gtrapprox K^{-1}\delta^{-v}.
\end{equation}

Consequently,
\begin{align*}
      I^C_{\mathbb{U}_4}(\mathbb{P}_4,\mathbb{B}_*) &\lesssim  WN|\Pcal_4|
     \lesssim_{\epsilon,C} W^{-3}N^{-1}\delta^{-O(\epsilon^2)}\delta^{-2}|\B_*|\\
     &\lesssim_{\epsilon,C} K^3\delta^{-O(\epsilon^2)}\delta^{2u-2}|\mathbb{B}_*|\\
     &\lesssim_{\epsilon,C} K^4 \delta^{1-\epsilon} |\mathbb{P}_4||\mathbb{B}_*|
\end{align*}
where the last inequality uses $2u+v\geq3$ .
. This proves the required estimate in the high--high case.

\medskip
\noindent\textbf{The low-frequency estimate \ref{low_2} holds.}

Since for each $U\in \M_3$,
\[
\#\{B\in\mathbb{B}_*:B\subseteq CU\}\sim N
\]
Thus by duality,
\begin{align}
    I^C_{\mathbb{U}_4}(\mathbb{P}_4,\mathbb{B}_*)  &\lesssim  I^C_{\mathbb{U}_3}(\mathbb{P}_4,\mathbb{B}_*) \sim N I^C(\M_3,\Pcal_4) \\
    &\lesssim N I^{D_0}(\hat{\Pcal}_4,\hat{\M}_3) \\
    &\lesssim_{\epsilon,C}
    S^{-1+\epsilon^2}N
    I_{\cap}^{C'}(\widehat{\mathbb U}_3,
                  \widehat{\mathbb P}_4^{S\delta}).
  \label{eq:dual-low-first-transfer}
    % &\lesssim_{\epsilon,C}    S^{-1+\epsilon^2} I_{\M_{\Delta}}^{100C^2_2C}(\Pcal_4^{S\delta},\mathbb{B}^{S\delta}_*)
\end{align}

To transfer the last count to $S\delta$ configuration, set
\begin{equation}\label{pair_2}
   \Phi:=\#\left\{ (p,B)\in \Pcal_4\times \B_* : \exists U\in\mathbb{U}_3,s.t.\; B\subseteq CU, \hat{p}_\Delta\cap   C'\widehat U\neq \emptyset \right\}
\end{equation}

On the one hand, 
for each $U\in\mathbb{ U}_3$ satisfies \ref{pair_2}, by construction of $\M'_\Delta$ in step 2, there exist $U_\Delta \in \M'_\Delta$, such that $U \subseteq C_* U_\Delta$. After duality, $\hat{U}\subseteq 10 C_* \widehat U_\Delta$. Since $\hat{p}_\Delta\cap C_*\hat{U}\neq \emptyset$, then
\[
\hat{p}_\Delta\subseteq\mathcal{N}_{10\Delta}(C_*\hat{U})\subseteq 100C_*\hat{U}_\Delta
\]
By Lemma \ref{dual_incidence_relation}, thus $U_\Delta \subseteq D_1 p^\Delta$, where $D_1=100C_*D_0$ is a fixed constant. Also, $B_\Delta \subseteq D_1U_\Delta$. Then
\begin{equation}\label{phi_l}
    \Phi\lesssim I_{\M_{\Delta}}^{D_1}(\Pcal_4^{S\delta},\mathbb{B}^{S\delta}_*)
\end{equation}
Here $\Pcal_4^{S\delta},\mathbb{B}^{S\delta}_*$ are counted with multiplicity.

On the other hand, we claim that if two boxes $U_1, U_2$ satisfy condition \eqref{pair_2} for the same pair $(p,B)$, then they must be $C^4$-comparable. Suppose otherwise. Let the directions of the long axes of $U_1$ and $U_2$ be $(\cos\theta_1, \sin\theta_1)$ and $(\cos\theta_2, \sin\theta_2)$, respectively. Since $B \subseteq CU_1 \cap CU_2$, if they are not comparable, then $|\theta_1 - \theta_2| \geq 10C\delta^{1/2}$. 

 Choose a point $(x_*, y_*) \in B$. By Definition \ref{definition_duality}, we have $ \hat{U_i} \subseteq C_2 V_{\theta_i}\cup C_2 V_{\theta_i+\pi}$, where $C_2$ is an absolute constant related to $C$  and
 \[
 V_\theta := \mathcal{N}_{\delta} \left\{ (x_* + z\sin\phi, y_* - z\cos\phi, z) : z \in [1,2], \, |\phi - \theta| \leq \delta^{1/2} \right\}.
 \]
Let $\delta\ll1$ such that $ \delta^{1/2}\gg C_2C'\delta
$. Consequently, $\mathcal{N}_{C'\Delta}(\hat{U}_1) \cap \mathcal{N}_{C'\Delta}(\hat{U}_2) = \emptyset$, which yields a contradiction. 
% Also see figure .

For each $U\in\M_3$, we have
\[
\#\{(p,B)\in \Pcal_4\times \B_* :  B\subseteq CU, \hat{p}_\Delta\cap C' \hat{U}\neq \emptyset \}\sim N |\{\hat{p}_\Delta\in \Pcal_4^{S\delta}:\hat{p}_\Delta\cap C' \hat{U}\neq \emptyset\}|
\]
Therefore, summing over all $U \in \M_3$, each pair $(p,B)$ is counted at most $C^{100}$ times. Thus,
\begin{equation}\label{Phi_g}
    \Phi \sim_{C} N
    I_{\cap}^{C'}(\widehat{\mathbb U}_3,
                  \widehat{\mathbb P}_4^{S\delta})
\end{equation}
Combining \ref{eq:dual-low-first-transfer} , \ref{phi_l} , \ref{Phi_g}, we have
\begin{equation}\label{eq:dual-low-contraction}
  I^C_{\mathbb U_4}(\mathbb P_4,\mathbb B_*)
  \lesssim_{\epsilon,C_0,C}
   S^{-1+\epsilon^2}I_{\M_{\Delta}}^{D_1}(\Pcal_4^{S\delta},\mathbb{B}^{S\delta}_*).
\end{equation}
Since $ \Pcal'_\Delta,\M'_\Delta,\B'_\Delta$  $(\Delta,C'_0,u,KK_0,i,j)$-nice configuration with $K_0\lessapprox1$,  by induction hypothesis on $\Delta=S\delta$,
\begin{equation}
    I^{D_1}_{\mathbb{U}'_\Delta}(\mathbb{P}'_{\Delta},\mathbb{B}_{\Delta,*}) \lesssim_{\epsilon,C_0,C} K_0^4K^4 (\delta S)^{1-\epsilon} |\mathbb{P}'_\Delta||\mathbb{B}_{\Delta,*}|
\end{equation}
Since $ \Pcal_1$ and $\B_*$ are $\{\delta,\Delta,1\}$ uniform, 
\begin{equation}
    I_{\M_{\Delta}}^{D_1}(\Pcal_4^{S\delta},\mathbb{B}^{S\delta}_*)\lesssim_{\epsilon,C_0,C} K_0^4K^4 (\delta S)^{1-\epsilon} |\mathbb{P}_1||\mathbb{B}_{*}|
\end{equation}
Combining with \ref{eq:dual-low-contraction},
\begin{equation}
    I^C_{\mathbb{U}_4}(\mathbb{P}_4,\mathbb{B}_*) \lesssim_{\epsilon,C_0,C} K_0^4K^4 S^{\epsilon^2-\epsilon
 }\delta ^{1-\epsilon} |\mathbb{P}_4||\mathbb{B}_{*}|
\end{equation}
Choosing $\delta<\delta_0\lesssim_{C,\epsilon} 1$, it is straightforward to close the induction.

\noindent\textbf{The low-frequency estimate \ref{low_1} holds.}

Since for each $ U\in \M_3$, 
\[
\#\{p \in \Pcal_2: U\subseteq Cp\}\sim M
\]
Thus
\begin{align}
    I^C_{\mathbb{U}_4}(\mathbb{P}_4,\mathbb{B}_*)  &\lesssim  I^C_{\mathbb{U}_3}(\mathbb{P}_2,\mathbb{B}_*) \sim M I^C(\M_3,\B_*) \\
    &\lesssim_{\epsilon,C} S^{-1+\epsilon^2} MI_{\cap}^{C'}(\mathbb U_3,\mathbb B_*^{S\delta}) 
    % &\lesssim_{\epsilon,C}    S^{-1+\epsilon^2} I_{\M_{\Delta}}^{C}(\Pcal_2^{S\delta},\mathbb{B}^{S\delta}_*)
\end{align}
Similarly, consider
\begin{equation}\label{pair_1}
   \Psi:=\#\{ (p,B)\in \Pcal_2\times \B_* : \exists U\in\mathbb{U}_3,\; U\subseteq Cp, B_\Delta\cap C' U\neq \emptyset \}
\end{equation}
On the one hand, 
for each $U\in\mathbb{ U}'_3$ satisfying \ref{pair_1}, we choose $U_\Delta \in \M'_\Delta$ such that $U \subseteq C_* U_\Delta$.  Since $B_\Delta\cap C'U\neq \emptyset$ and $U\subseteq Cp$, there exits a constant $D_2$ (depending on $C,C',C_*$) such that
\[
B_\Delta\subseteq D_2  U_\Delta,\qquad U_\Delta \subseteq D_2 p^\Delta
\]
Hence,
\begin{equation}
    \label{Psi_l}
    \Psi \lesssim I_{\M_{\Delta}}^{D_2}(\Pcal_2^{S\delta},\mathbb{B}^{S\delta}_*)
\end{equation}

On the other hand, similarly, if two boxes $U_1,U_2$ satisfies \ref{pair_1} for the same $(p,B)$, they must be $C^4$-comparable. Suppose that $U_1, U_2 \subseteq Cp$ are distinct. Then they are contained in $\delta^{1/2}$-arcs $\omega_1, \omega_2$ on the $\delta$-circle $p$, which are separated by at least $C\delta^{1/2}$. It follows that
\[
B_\Delta \subseteq \mathcal{N}_{C'\Delta}(\omega_1) \cap \mathcal{N}_{C'\Delta}(\omega_2),
\]
which yields a contradiction.  
% Also see figure .

 For $U\in\M_3$, 
\[
\#\{(p,B)\in \Pcal_2\times \B_* : U\subseteq Cp, B_\Delta\cap C'U\neq \emptyset \}\sim M |\{B\in \B_\Delta:B_\Delta\cap C' U\neq \emptyset\}|
\]
So summing up for $U\in\M_3$, we have
\begin{equation}
    \label{Psi_g}
   \Psi \sim_C MI_{\cap}^{C'}(\mathbb U_3,\mathbb B_*^{S\delta})
\end{equation}
Combining the preceding estimates,
\begin{equation}
  I^C_{\mathbb U_4}(\mathbb P_2,\mathbb B_*)
  \lesssim_{\epsilon,C_0,C}
   S^{-1+\epsilon^2}I_{\M_{\Delta}}^{D_2}(\Pcal_2^{S\delta},\mathbb{B}^{S\delta}_*).
\end{equation}
Similarly by induction, we have
\[
I_{\M_{\Delta}}^{D_2}(\Pcal_2^{S\delta},\mathbb{B}^{S\delta}_*)\lesssim_{\epsilon,C_0,C} K_0^4K^4 (\delta S)^{1-\epsilon} |\mathbb{P}_1||\mathbb{B}_{*}|
\]
And we can close the induction as above.

\section{proof for the cinematic case}
\label{sec:cinematic-extension}

In this section, we prove Theorem \ref{Thm_discretized} for the cinematic curves defined in \ref{def:cinematic-family}. The frame of the proof remains the same with the circular case. The new ingredient is a variable coefficient high low argument. 

We work in one compact chart from Definition~\ref{def:cinematic-family},
with fixed buffers $Z\Subset Z_1\Subset Z^+$ and
$I\Subset I_1\Subset I^+$. We usually use the notation $q\in \mathbb Q$ representing $\delta$ discretized ball in parameter space $Z$ and $p_q\in\Pcal$ referring its corresponding $\delta$ discretized cinematic curve.

Similarly with arguments in Section ~5, we introduce the dual incidence problem between tangent boxes and cinematic curves.  

\subsection{Duality for the cinematic case}

For a tangent box, use the graph-coordinate form
\[
 U=U_\delta(t_U,y_U,s_U)
 :=\{(t,y):|t-t_U|\leq\delta^{1/2},\ 
       |y-y_U-s_U(t-t_U)|\leq\delta\}.
\]
The following is the cinematic analogue of
Definition~\ref{definition_duality}. Since we consider pairwise distinct tangent boxes in a bounded region, we can assume that $t_U,s_U\in \delta^{1/2}\mathbb Z \cap [-1,1], y_U \in \delta\mathbb Z \cap [-1,1]$.

\begin{definition}[Cinematic duality plank]
\label{def:cinematic-duality}
Define
\begin{equation}\label{eq:cinematic-duality-plank}
 \Gamma_U^\delta
 :=\{z\in Z_1:|g(z,t_U)-y_U|\leq\delta,
                 \ |g_t(z,t_U)-s_U|\leq\delta^{1/2}\}.
\end{equation}
and
\begin{equation}\label{eq:cinematic-duality-plank_variant}
 \Gamma_U^{\alpha,\beta}
 :=\{z\in Z_1:|g(z,t_U)-y_U|\leq\alpha,
                 \ |g_t(z,t_U)-s_U|\leq\beta\}.
\end{equation}
For $C>0$, $C\Gamma_U^\delta$ denotes the enlargement set $\Gamma_U^{C\delta,C\delta^{1/2}}$. 

% Its central locus
% $\{g(z,t_U)=y_U,\ g_t(z,t_U)=s_U\}$ may be curved, so this notation
% does not mean Euclidean dilation.
\end{definition}

The following is the cinematic version of Lemma \ref{lem:duality}. 
\begin{lemma}[Cinematic dual incidence relation]
\label{lem:cinematic-duality}
For each fixed $C\geq1$ there is $C_1=C_1(C,g)$ such that, for a
parameter cube $q$, its discretized curve $p_q$, and a tangent box $U$,
\begin{equation}\label{eq:cinematic-dual-incidence}
 U\subset Cp_q\ \Longrightarrow\ q\subset C_1\Gamma_U^\delta,
 \qquad
 q\cap C\Gamma_U^\delta\ne\varnothing\ \Longrightarrow\
 U\subset C_1p_q.
\end{equation}
Moreover, for fixed $C$ and $\delta\leq r\leq1$,
\begin{equation}\label{eq:cinematic-dual-volume}
 |C\Gamma_U^\delta|\lesssim_{C,g}\delta^{3/2},\qquad
 |C\Gamma_U^\delta\cap B(z,r)|\lesssim_{C,g}\delta r^2.
\end{equation}
For $\delta<\Delta\ll\delta^{1/2}$,
\begin{equation}\label{eq:enlarge_plank}
    \Gamma^{\Delta,\delta^{1/2}}_U \subseteq \mathcal{N}_{O(\Delta)} (\Gamma_U^\delta)
\end{equation}
\end{lemma}
\begin{proof}

Fix \(z\in q\). Suppose first that \(U\subset Cp_q\). By Lemma \ref{lem:cinematic-jet-geometry}, $U\subset \mathcal{N}_{CC_g\delta}(\Sigma_z)$. For
\(j\in\{0,1\}\), let
\[
 t_j=t_U+j\delta^{1/2},\qquad Y_j=y_U+js_U
 \delta^{1/2}.
\]
 \((t_j,Y_j)\in U\), for each \(j\) give that 
 \[
 |g(z,t_U)-y_U|\leq CC_g \delta, \quad |g(z,t_1)-g(z,t_{0})-s_U\delta^{1/2}|\leq CC_g \delta
 \]
 Since $|g(z,t_1)-g(z,t_{0})-g_t(z,t_U)(t_1-t_0)|\lesssim \delta$. Thus
\[
 |g(z,t_U)-y_U|\leq CC_g \delta, \quad |g_t(z,t_U)-s_U|\lesssim_g C \delta^{1/2}
\]
 Thus \(z\in C_1\Gamma_U^\delta\) for some $C_1=C_1(C,g)$. Conversely, it is similar.

By lemma \ref{lem:cinematic-jet-geometry}, \[
  J_t:z\longmapsto\bigl(g(z,t),g_t(z,t),g_{tt}(z,t)\bigr)
\]
are uniformly bi-Lipschitz. Since $J_t(C\Gamma_U^\delta)$ is contained in a box with dimension $C\delta\times C\delta^{1/2}\times C_g$, and  $J_t[C\Gamma_U^\delta\cap B(z,r)]$ is contained in a box with dimension $ C\delta\times C_g r\times C_g r$, where $C_g =\|g\|_{C^3(Z\times I)}$. Thus \ref{eq:cinematic-dual-volume} holds.

$\Gamma_U^\delta=J^{-1}_t(V)\cap Z$ for a rectangle $V$ with dimension $\delta\times \delta^{1/2}$. And $\Gamma_U^{\Delta,\delta^{1/2}}$ is contained in a $J^{-1}_t[\mathcal{N}_{2\Delta}(V)]$. By the bi-Lipschitz property,  which is in $\mathcal{N}_{O(\Delta)} (\Gamma_U^\delta)$.
\end{proof}

\subsection{The cinematic high--low method}

We retain the notation of Section~3 for the cinematic phase
\eqref{eq:cinematic-phase}, the annular wedge $\mathcal A$, the Fourier
integral operator \eqref{eq:cinematic-fio}, and its wave packets. The main point is the passage from high frequency part of characteristic function for dual planks to a common Fourier integral operator.  In particular, each dual plank $\Gamma^\delta_U$ is approximately a wave packet of $T_\lambda$ after rescaling.

% All Fourier transforms use
% the convention
% $\widehat h(\xi)=\int e^{-2\pi i x\cdot\xi}h(x)\,dx$.

Fix $0<\epsilon\ll1$, and put
\begin{equation}\label{eq:cin-hl-parameters}
 R=\delta^{-1},\qquad S=R^\epsilon,\qquad
 \kappa=\epsilon/4,\qquad S_1=S^{1-\kappa}.
\end{equation}
Let $A=R^\rho$, where $\rho>0$ is chosen sufficiently small below;
in particular, $A^{10}\ll S^\kappa$.

For $U=U_\delta(t_U,y_U,s_U)$, write
\[
 v_U(z)=g(z,t_U)-y_U,\qquad
 w_U(z)=g_t(z,t_U)-s_U.
\]

Let $\chi$ be a smooth chart cutoff on $Z$, and nonnegative $\chi_0,\chi_1\in C^\infty_c(\R)$. Define
\begin{equation}\label{eq-cin-fU}
 f_U(z)=\chi(z)\chi_0(v_U(z)/\delta)
                    \chi_1(w_U(z)/\sqrt\delta).
\end{equation}
In particular, $f_U\geq1$ on the relevant portion of
$C\Gamma_U^\delta$. Fourier inversion gives
\begin{equation}\label{eq:cin-fU-inversion}
 \begin{split}
 f_U(z)
 ={}&\delta^{3/2}\chi(z)\iint
 \widehat\chi_0(\delta\sigma)
 \widehat\chi_1(\sqrt\delta\,\tau)\\
 &\hspace{12mm}\times
 e^{2\pi i[\sigma v_U(z)+\tau w_U(z)]}\,d\sigma\,d\tau .
 \end{split}
\end{equation}
Choose a bump function \(\varphi \in C_c^\infty(\mathbb{R})\) which equals 1 on $B_1$ and vanishes out of $B_2$. Define low frequency part \(f_U^\ell\) by inserting \(\varphi(S_1\delta\sigma)\) in \ref{eq:cin-fU-inversion}, and set \(f_U^h = f_U - f_U^\ell\). Decompose \(1 - \varphi(S_1\delta\sigma)\) into smooth dyadic shell multipliers \(\psi_\lambda(\sigma)\), supported where \(|\sigma| \sim \lambda\). We write \(f_U^\lambda\) for the corresponding term. Since $\widehat\chi_0$ and $\widehat\chi_1$ are
Schwartz, we may restrict to
\begin{equation}\label{eq:cin-lambda-range}
    c(S_1\delta)^{-1} \le \lambda \le CA\delta^{-1} 
\end{equation}
In particular, 
\begin{equation}
    f_U=f^l_U+f^h_U=f^l_U+\sum_{ \lambda}f^\lambda_U+RapDec(R)
\end{equation}

\begin{proposition}[Transfer to FIO with cinematic phase]
\label{prop:cin-fio-transfer}
Fix $\lambda$ in \eqref{eq:cin-lambda-range}, and a parameter ball
$B=B(z_B,c/S)$. Put $ r=\lambda /S$. For $z\in2B$, we have
\begin{equation}\label{eq:cin-transfer}
 f_U^\lambda(z)
 =\lambda^2\delta^{3/2}
 \sum_{|m|\leq A}e^{2\pi i m\cdot S(z-z_B)}
 T_{\lambda,B}h_{U,m}\bigl(\lambda(z-z_B)\bigr)
 {}+RapDec(R),
\end{equation}
where the operator is defined as 
\begin{equation}\label{eq:cin-local-fio}
 \begin{split}
 T_{\lambda,B}h(y)
 &=\int e^{2\pi i\lambda[
       \phi(z_B+y/\lambda,\xi)-\phi(z_B,\xi)]}
       a(z_B+y/\lambda,\xi)h(\xi)\,d\xi .
 \end{split}
\end{equation}
Here $a(z,\xi)=\chi(z)\eta(\xi)$, where $\eta$ is a fixed
annular-wedge cutoff on $\mathcal{A}$.

Set
\[
 y_{U,z_B}=y_U-g(z_B,t_U),\qquad s_{U,z_B}=s_U-g_t(z_B,t_U).
\]
And
\begin{equation}\label{eq:cin-hUm}
 \begin{split}
 h_{U,m}(\xi)
 ={}&c_{m,U}(\xi)\psi_\lambda(\lambda\xi_2)
       \widehat\chi_0(\delta\lambda\xi_2)\\
 &\times\widehat\chi_1\bigl(
        \sqrt\delta\,\lambda(\xi_1-t_U\xi_2)\bigr)\\
 &\times e^{-2\pi i\lambda[
        \xi_2 y_{U,z_B}+(\xi_1-t_U\xi_2) s_{U,z_B}]}.
 \end{split}
\end{equation}
Thus, $h_{U,m}$ is essentially supported on $|\xi_2|\sim1$, $|\xi_1/\xi_2-t_U|\lesssim r^{-1/2}$. Writing $t=\xi_1/\xi_2$, for any $N>1$, the coefficients $c_{m,U}$ satisfy
\begin{equation}\label{eq:cin-cm-symbol}
 \left|\partial_{\xi_2}^j
       (r^{-1/2}\partial_t)^l c_{m,U}\right|
 \lesssim_{N,j,l,g}(1+|m|)^{-N}.
\end{equation}
\end{proposition}
\begin{proof}

\smallskip\noindent
\emph{Local normalization.}
Fix a parameter ball \(B(z_B,c/S)\). By performing the following local normalization, we may reduce to the case where \(z_B=0\) and \(\phi(0,\cdot)\equiv 0\). Let \(w=z-z_B\) and set
\[
 g^{z_B}(w,t)=g(z_B+w,t)-g(z_B,t),\qquad
 \phi^{z_B}(w,\xi)=\phi(z_B+w,\xi)-\phi(z_B,\xi).
\]
Then \(g^{z_B}(0,t)\equiv0\) and \(\phi^{z_B}(0,\xi)\equiv0\). Furthermore, let
\(
 a^{z_B}(w,\xi)=a(z_B+w,\xi),
\)
so that \(T^{z_B}_{\lambda,B} h=T_{\lambda,B} h\). Hence, in the proof, we assume \(z_B=0\) and \(\phi(0,\cdot)\equiv 0\).

Make the determinant-one change of variables \(\zeta_2=\sigma\), \(\zeta_1=t_U\sigma+\tau\), and put \(t=\zeta_1/\zeta_2\). Define
\[
 q_U(z,\theta)=g(z,t_U)+(t-t_U)g_t(z,t_U)-g(z,t).
\]
Comparing the exact phase in \ref{eq:cin-fU-inversion} with the cinematic phase \(\phi\), we obtain
\begin{equation}\label{eq:cin-phase-identity}
 \begin{split}
 \sigma v_U(z)+\tau w_U(z)
 ={}&\phi(z,\zeta)-\sigma y_U-\tau s_U
       +\sigma q_U(z,t).
 \end{split}
\end{equation}
After the substitution \(\zeta=\lambda\xi\), the last term becomes
\[
 E_U(z,\xi)
 =\lambda\xi_2q_U(z,\theta),\quad t=\xi_1/\xi_2.
\]
The normalization \(g(0,t)\equiv0\) implies
\(\partial_t^j g(z,t)=O_g(|z|)\) for each fixed \(j\). Taylor's theorem
therefore yields \(q_U=O_g(|z|\,|t-t_U|^2)\). Since
\(|z|\lesssim S^{-1}\),
\[
 |E_U|\lesssim_g
 \lambda S^{-1}r^{-1}\lesssim_g1
 \quad\text{when }|\theta-t_U|\lesssim r^{-1/2}.
\]
The same estimate holds after applying any fixed number of derivatives in
\(u=Sz\), \(\xi_2\), and \(r^{1/2}(t-t_U)\).

Next, expand \(e^{2\pi iE_U}\) in a Fourier series in \(u=Sz\in B_1\). On \(2B\), this gives
\[
 e^{2\pi iE_U(z,\xi)}
 =\sum_{m\in\mathbb Z^3}
       e^{2\pi im\cdot Sz}c_{m,U}(\xi).
\]
Integration by parts proves \eqref{eq:cin-cm-symbol}. Consequently, it follows that
\begin{equation}\label{eq:cin-fU}
 \begin{split}
 f^\lambda_U(z)
 ={}&\lambda^2\delta^{3/2}\chi(z)\sum_{m\in\mathbb Z^3}\iint
c_{m,U}(\xi)\psi_\lambda(\lambda\xi_2)
       \widehat\chi_0(\delta\lambda\xi_2)
\\
&\times\widehat\chi_1\bigl(
        \sqrt\delta\,\lambda(\xi_1-t_U\xi_2)\bigr) e^{2\pi i\lambda\phi(z,\xi)-2\pi i\lambda[
        \xi_2 y_{U}+(\xi_1-t_U\xi_2) s_{U}]} d\xi,
 \end{split}
\end{equation}
where the factor \(\lambda^2\) is the Jacobian of \(\zeta=\lambda\xi\).

Finally, after absorbing the terms with \(|m|>A\) into the \(RapDec\) terms, \ref{eq:cin-transfer} follows.
    
\end{proof}

The following proposition shows that \(h_{U,m}\) is approximately a single wave packet associated with \(T_U\) for \(T_{\lambda,B}\). Moreover, after a suitable translation and dilation, \(T_U\) is exactly the dual plank \(\Gamma^\delta_U\).

\begin{proposition}[wave packets and dual planks]
\label{prop:cin-packet-source}
Use the standard $r$ scale wave packet decomposition from
Section ~3 for $T_{\lambda,B}$
on $|y|\leq r$. After discarding negligible
tails, each $h_{U,m}$ has a collection of wave packets $\mathbb T_{\lambda,B}[U]$,
 chosen uniformly for $|m|\leq A$, with
the following properties:
\begin{enumerate}[label=\textup{(\roman*)}]
 \item $|\mathbb T_{\lambda,B}[U]|\leq S^{O(1)}$. And each wave packet $T\in\mathbb T$ belongs to $ \mathbb T_{\lambda,B}[U]$ for at most $S^{O(1)}$ distinct tangent boxes $U\in\M$.
 \item If $T\in\mathbb T_{\lambda,B}[U]$ and
 $z\in 2B$ with $\lambda(z-z_B)\in T$, then
 \begin{equation}\label{eq:cin-source-jets}
 |v_U(z)|\lesssim_g A\lambda^{-1},\qquad
 |w_U(z)|\lesssim_g\sqrt\delta.
 \end{equation}
 In particular, 
 \begin{equation}\label{eq:cin-source-neighborhood}
 \bigl(z_B+\lambda^{-1}T\bigr)\cap 2B
 \subset\mathcal N_{o(S\delta)}(C'\Gamma_U^\delta)
 \end{equation}
 for a fixed sufficiently large $C'=C'(g,C)$.

\end{enumerate}
\end{proposition}

\begin{proof}

Similarly, we assume \(z_B=0\) and \(\phi(0,\cdot)\equiv 0\). Choose \(\epsilon_{\rm dec}\ll\rho\ll\epsilon\) and set
\begin{equation}\label{eq:cin-dec-scales}
 1\ll r=\lambda/S\leq\lambda^{1-\epsilon_{\rm dec}},
 \qquad
 r^{\epsilon_{\rm dec}},\ r^{2\nu},\ K_{\rm wp}\leq A,
\end{equation}
where \(\nu=\epsilon_{\rm dec}^{1000}\), \(K_{\rm wp}=r^{\epsilon_{\rm dec}^{50}}\), and \(A=R^\rho\).

Denote by \(\mathbb T=\{T_{\theta,V}\}\) the \(r\)-scale wave packets in \(B_r\). Recalling Section 3,
\begin{equation}
 \theta=\left\{\xi\in\mathcal A:
 \left|\frac{\xi_1}{\xi_2}-t_\theta\right|\leq r^{-1/2}\right\},
 \qquad
 t_\theta\in r^{-1/2}\mathbb Z.
\end{equation}
For each \(V\in\mathcal V_\theta\), \(s_V\in\lambda^{-1}r^{1/2}\mathbb Z\), and \(y_V\in\lambda^{-1}\mathbb Z\),
\begin{equation}
 \lambda^{-1}V=\{(s,y):|s-s_V|\leq\lambda^{-1}r^{1/2},\ |y-y_V+t_\theta s|\leq\lambda^{-1}\}.
\end{equation}

For \(h_{U,m}\), Proposition \ref{prop:cin-fio-transfer} and \eqref{eq:cin-hUm} show that its frequency part is essentially supported on
\begin{equation}\label{eq:frequency-area}
 |\xi_2|\sim1,\quad |\xi_1/\xi_2-t_U|\lesssim r^{-1/2},
\end{equation}
and its physical part \(\check h_{U,m}\) is essentially supported on \(\lambda\Lambda_U\), where \(\Lambda_U\) is the rectangle
\begin{equation}
 \{(s,y):|s-s_U|\leq\delta^{1/2},\ |y-y_U+t_Us|\leq\lambda^{-1}\}.
\end{equation}
Thus the essential wave packets of \(h_{U,m}\) are \(\mathbb T_{\lambda,B}[U]=\{T_{\theta,V}\}\), where \(\theta\) intersects \eqref{eq:frequency-area} and \(\tilde V\) intersects \(\lambda\Lambda_U\). It follows that \((\theta,V)\) satisfies
\begin{equation}\label{eq:approx}
 |t_\theta-t_U|\lesssim r^{-1/2},\qquad |s_U-s_V|\leq\delta^{1/2},
\end{equation}
since \(r^{\nu}\lambda^{-1}r^{1/2}\leq\delta^{1/2}\), and
\begin{equation}
 |y_U-y_V-(t_U-t_\theta)s_V|\lesssim r^{2\nu}\lambda^{-1}.
\end{equation}
Hence \(\#(\theta,V)\lesssim S^{O(1)}\). Also, fix \(T=T_{\theta,V}\); for a distinct family \(\M\) with \(t_U,s_U\in\delta^{1/2}\mathbb Z\) and \(y_U\in\delta\mathbb Z\), there are at most \(S^{O(1)}\) choices of \(U\in\M\) satisfying the above. Thus (i) holds.

Recalling \ref{eq_shape_of_wavepacket},
\begin{equation}
 \lambda^{-1}T_{\theta,V}=\{z\in B:|g_t(z,t_\theta)-s_V|\leq r^{1/2+\nu}/\lambda,\ |g(z,t_\theta)-y_V|\leq r^{2\nu}/\lambda\}.
\end{equation}
Then for \(z\in\lambda^{-1}T\),
\begin{align*}
 |g_t(z,t_U)-s_U-g_t(z,t_\theta)+s_V|
 &\leq |g_t(z,t_\theta)-g_t(z,t_U)|+|s_U-s_V|\\
 &\lesssim S^{-1}|t_\theta-t_U|+\delta^{1/2}\lesssim\delta^{1/2},
\end{align*}
and
\begin{align*}
 |g(z,t_U)-y_U-g(z,t_\theta)+y_V|
 &\leq |g(z,t_U)-g(z,t_\theta)-g_t(z,t_\theta)(t_U-t_\theta)|\\
 & +|g_t(z,t_\theta)-s_V||t_U-t_\theta|+|y_U-y_V-(t_U-t_\theta)s_V|\\
 &\lesssim S^{-1}r^{-1}+r^{\nu}/\lambda+r^{2\nu}/\lambda\lesssim A/\lambda.
\end{align*}
Thus \ref{eq:cin-source-jets} follows. Equivalently,
\[
 \lambda^{-1}T\subseteq C\Gamma_U^{A\lambda^{-1},\delta^{1/2}}.
\]
Since \(A\lambda^{-1}\lesssim AS_1\delta=o(S\delta)\), Lemma \ref{lem:cinematic-duality} implies \ref{eq:cin-source-neighborhood}.

\end{proof}

Next, we prove the cinematic high-low estimate (cinematic version of Theorem \ref{thm_high_low_2}) to bound the incidence between $\delta$ discretized ball in $\R^3$ and curved planks $\{\Gamma^\delta_U\}$.

\begin{theorem}[Cinematic high--low estimate]\label{thm:cinematic-high-low}
Let $\mathbb Q$ be a family of distinct dyadic parameter $\delta$-cubes
in the fixed chart, and let $\mathbb U$ be a distinct family of tangent
boxes. Fix $C\geq1$, and suppose that each
$q\in\mathbb Q$ is contained in $\sim W$ planks $C\Gamma_U^\delta$.
Put
\begin{equation}\label{eq:cin-IJ}
 \begin{split}
 \mathcal I&=\sum_{q\in\mathbb Q}
       \#\{U\in\mathbb U:q\subset C\Gamma_U^\delta\},\\
 \mathcal J&=\sum_{q\in\mathbb Q}
       \#\{U\in\mathbb U:
            q_{S\delta}\cap C'\Gamma_U^\delta
                     \ne\varnothing\},
 \end{split}
\end{equation}
where $C'=C'(C,g)$ is fixed and sufficiently large.
Then at least one of the following holds: 
\begin{enumerate}
    \item The high-frequency dominates
\begin{equation} \label{eq:cin-high-alternative}
    \mathcal I
 \lesssim_{\epsilon,C,g}
 S^{O(1)}W^{1/3}\delta^{-1/4}
             |\mathbb Q|^{5/6}|\mathbb U|^{1/6},
\end{equation}
Equivalently,
\begin{equation}\label{eq:cin-richness-alternative}
 W^4|\mathbb Q|
       \lesssim_{\epsilon,C,g}S^{O(1)}\delta^{-3/2}|\mathbb U|.
\end{equation}
\item The low-frequency dominates
\begin{equation} \label{eq:cin-low-alternative}
    \mathcal I
 \lesssim_{\epsilon,C,g}S^{-1+\epsilon} \mathcal  J
\end{equation}

\end{enumerate}
\end{theorem}

\begin{proof}
Put \(f=\sum_U f_U\) and, for \(\mathcal A\subseteq\mathbb Q\), put
\(G_{\mathcal A}=\sum_{q\in\mathcal A}\one_q\). Let
\[
 \mathcal I_{\mathcal A}=\sum_{q\in\mathcal A}
 \#\{U\in\mathbb U:q\subset C\Gamma_U^\delta\}.
\]
By \ref{eq-cin-fU},
\begin{equation}\label{eq:cin-incidence-pairing}
 \mathcal I_{\mathcal A}\sim W|\mathcal A|\lesssim
 \delta^{-3}\int fG_{\mathcal A}.
\end{equation}

Write \(f=f^\ell+f^h\), where \(f^\ell=\sum_U f^\ell_U\). We have
\[
 f_U^\ell(z)
 =\chi(z)b_\ell(v_U(z))
 \chi_1(w_U(z)/\sqrt\delta),
\]
where
\[
 b_\ell(v)=\delta\int
 \widehat\chi_0(\delta\sigma)\varphi(S_1\delta\sigma)
 e^{2\pi i\sigma v}\,d\sigma.
\]
The function \(b_\ell\) is essentially supported on an \(S_1\delta\) ball and is equal to \(S_1^{-1}\) there. More precisely, choosing \(S_2=S^{1-\kappa/2}\in(S_1,S)\),
\begin{equation}\label{eq:cin-low-profile}
 |b_\ell(v)|\lesssim_N
 S_1^{-1}\left(1+\frac{|v|}{S_1\delta}\right)^{-N},\quad |v|>S_2\delta,\qquad |b_\ell(v)|=\mathrm{RapDec}(R).
\end{equation}
Therefore \(f^\ell_U\) is essentially supported on \(C\Gamma^{S_2\delta,\delta^{1/2}}_U\), by Lemma \ref{lem:cinematic-jet-geometry}, which is contained in \(\mathcal N_{O(S_2\delta)}(C'\Gamma_U)\subseteq \mathcal N_{S\delta}(C'\Gamma_U)\). Denote
\[
 m_S(q):=\#\{U:q_{S\delta}
 \cap C'\Gamma_U^\delta\ne\varnothing\}.
\]
Consequently, for every \(\mathcal A\),
\begin{equation}\label{eq:cin-low-pairing}
 \delta^{-3}\left|\int f^\ell G_{\mathcal A}\right|
 \lesssim_{\epsilon,C,g}
 S^{-1+\kappa}\sum_{q\in\mathcal A}m_S(q)+\mathrm{RapDec}(R).
\end{equation}

If the low estimate \eqref{eq:cin-low-alternative} fails with a sufficiently large constant, then
\[
 \mathcal J\sim \sum_{q\in\mathbb Q} m_S(q) \lesssim S^{1-\epsilon}W|\mathbb Q|.
\]
Discard cubes with \(m_S(q)>2\mathcal{J}/|\mathbb Q|\). The remaining set \(\mathbb Q_0\) has at least \(|\mathbb Q|/2\) elements and satisfies
\begin{equation}\label{eq:cin-good-cubes}
 \sup_{q\in\mathbb Q_0}m_S(q)
 \lesssim_{\epsilon,C,g}S^{1-\epsilon}W.
\end{equation}
Also, by \eqref{eq:cin-low-pairing}, the high frequency part dominates, and hence
\begin{align}
 \mathcal I &\sim \mathcal I_{\mathbb Q_0} \lesssim_{\epsilon,C,g} \delta^{-3}\left|\int f^h G_{\mathbb Q_0}\right|\\
 &\leq \delta^{-3}\|f^h\|_{L^6(X)}\|G_{\mathbb Q_0}\|_{6/5}\label{eq_holder-high}
\end{align}
with \(X=\bigcup_{q\in\mathbb Q_0}q\).

Decompose \(f^h\) into the \(O(\log R)\) essential shells in \eqref{eq:cin-lambda-range}, and cover the chart \(Z\) by \(O_g(S^3)\) balls \(B\) of radius \(c/S\). Then
\begin{equation}
 \|f^h\|_{L^6(X)}\lesssim \sum_{\lambda,B} \|f^\lambda\|_{L^6(X\cap B)}+\mathrm{RapDec}(R).
\end{equation}

Fix \(\lambda\in (c(S_1\delta)^{-1},CA\delta^{-1})\) and \(B=B(z_B,c/S)\). Retain the symbols above and set
\begin{equation}
 1\ll r=\lambda/S\leq\lambda^{1-\epsilon_{\rm dec}},
 \qquad
 r^{\epsilon_{\rm dec}},\ r^{2\nu},\
 K_{\rm wp}\leq A,\qquad A^{10}\ll S^\kappa,
\end{equation}
where \(\nu=\epsilon_{\rm dec}^{1000}\), \(K_{\rm wp}=r^{\epsilon_{\rm dec}^{50}}\), and \(A=R^\rho\).

Apply Proposition~\ref{prop:cin-fio-transfer} on each ball. For each \(U\in\M\) and \(z\in B\),
\begin{equation}
 f_U^\lambda(z)
 =\lambda^2\delta^{3/2}
 \sum_{|m|\leq A}e^{2\pi i m\cdot S(z-z_B)}
 T_{\lambda,B}h_{U,m}\bigl(\lambda(z-z_B)\bigr)
 +\mathrm{RapDec}(R).
\end{equation}
Hence,
\begin{equation}
 \|f^\lambda\|_{L^6(X\cap B)} \lesssim \lambda^{3/2}\delta^{3/2} \sum_{|m|\leq A}\left\|\sum_{U\in\M}T_{\lambda,B} h_{U,m}\right\|_{L^6(X^r_B)}+\mathrm{RapDec}(R),
\end{equation}
where \(X^r_B:=\lambda(X\cap B-z_B)\) is a union of \(\lambda\delta\) balls in \(B_r\).

We perform the \(r\)-scale wave packet decomposition for \(\sum_{U\in\M}T_{\lambda,B} h_{U,m}\). Denote by \(\mathbb T\) the set of wave packets and by \(h_T\), \(T\in\mathbb T\), the corresponding functions. Recalling Proposition~\ref{prop:cin-packet-source}, after discarding negligible tails we have
\[
 h_T=\sum_{U:T\in\mathbb T_{\lambda,B}[U]}(h_{U,m})_T.
\]

Cover \(X^r_B\) by \(K_{\rm wp}\)-balls, denoted by \(\bar X_B^r\). Let \(L_B\) be the linear transform from \(B_r\) to \(B\). Then \(L_B(\bar X_B^r)\) contains a union of \(K_{\rm wp}/\lambda=o(S\delta)\) enlargements of balls from \(X\cap B\). Moreover, by Proposition \ref{prop:cin-packet-source}, \(|\mathbb T_{\lambda,B}[U]|\leq S^{O(1)}\), and for each \(T\in \mathbb T_{\lambda,B}[U]\),
\[
 L_B(T) \cap 2B
 \subset\mathcal N_{o(S\delta)}(C'\Gamma_U^\delta).
\]
Combining this with the multiplicity bounds \eqref{eq:cin-good-cubes}, each \(K_{\rm wp}\)-radius ball in \(\bar X_B^r\) intersects at most \(S^{O(1)}W\) wave packets. After discarding negligible wave packets \(T\) such that \(\|h_T\|_{L^2}\leq R^{-1000}\), we perform pigeonholing and assume that each \(\|T_{\lambda,B}h_T\|_{L^6(\omega_{B_r})}\) is the same up to a constant. Applying the variable coefficient refined decoupling for \(T_{\lambda,B}\),
\begin{equation*}
 \left\|\sum_{U\in\M}T_{\lambda,B} h_{U,m}\right\|_{L^6(\bar X^r_B)} \lesssim_{\epsilon} r^{\epsilon_{\rm dec}} W^{1/3} \left(
 \sum_{T\in\mathbb T}
 \|T_{\lambda,B} h_T\|_{L^6(\omega_{B_r})}^{6}
 \right)^{1/6}
 +\operatorname{RapDec}(r)\|h\|_2.
\end{equation*}
By Proposition \ref{prop:cin-packet-source} (i) and ~\ref{prop_wave_packet} (3),
\begin{align*}
 \|T_{\lambda,B} h_T\|^6_{L^6(\omega_{B_r})} &\leq S^{O(1)} \sum_{U:T\in\mathbb T_{\lambda,B}[U]} \|T_{\lambda,B}(h_{U,m})_T\|^6_{L^6(\omega_{B_r})}\\
 &\lesssim S^{O(1)} \sum_{U:T\in\mathbb T_{\lambda,B}[U]} \|h_{U,m}\|^6_{L^2},
\end{align*}
and
\begin{align*}
 \sum_{T\in\mathbb T}
 \|T_{\lambda,B} h_T\|_{L^6(\omega_{B_r})}^{6} \lesssim S^{O(1)} \sum_{U\in\M} \|h_{U,m}\|^6_{L^2}.
\end{align*}
Recalling the representation of \(h_{U,m}\) in Proposition \ref{prop:cin-fio-transfer}, we have
\[
\|h_{U,m}\|_{L^2}\lesssim \lambda^{-1/2}\delta^{-1/4}
\]
Combining the above, we obtain
\[
 \|f^\lambda\|_{L^6(X\cap B)} \lesssim_{\epsilon,C,g}S^{O(1)}W^{1/3}\delta^{1/4}|\mathbb U|^{1/6}
 +\mathrm{RapDec}(R).
\]
Since \(\#\{\lambda,B\}\lesssim S^4\), we have
\[
 \|f^h\|_{L^6(X)} \lesssim_{\epsilon,C,g}S^{O(1)}W^{1/3}\delta^{1/4}|\mathbb U|^{1/6}.
\]
By \eqref{eq_holder-high} and \(\|G_{\mathbb Q_0}\|_{6/5}\leq |\mathbb Q|^{5/6}\delta^{5/2}\), the high alternative \eqref{eq:cin-high-alternative} holds.
\end{proof}

\subsection{The configuration estimate}

\begin{theorem}[Cinematic configuration estimate]\label{Thm_configuration_cinematic}
Let $0<\epsilon\ll1$, $C_0,C>10$, $0<u\leq1$, and $0<v\leq3$, with
$2u+v\geq3$. There is $\delta_0>0$ such that the following holds for every
dyadic $\delta\in(0,\delta_0)$ and $K\geq1$. Suppose
$(\mathbb P_0,\mathbb U_0,\mathbb B_0)_\delta$ is a 
$(\delta,C_0,u,K,I,J)$-nice configuration of cinematic graphs whose parameter labels forms a
$(\delta,v,K)$-set. Then there are a refinement
$(\mathbb P',\mathbb U',\mathbb B')_\delta$ and a family
$\mathbb B'\subseteq\mathbb B_*\subseteq\mathbb B_0$,
\begin{equation}\label{eq:configuration-incidence-bound-cin}
  I^C_{\mathbb U'}(\mathbb P',\mathbb B_*)
  \lesssim_{\epsilon,C_0,C,g}
  K^4\delta^{1-\epsilon}|\mathbb P'||\mathbb B_*|.
\end{equation}
\end{theorem}

\begin{proof}
    We follow the four steps of the proof of Theorem \ref{Thm_configuration}, and provide the corresponding replacements in the cinematic setting wherever the circle geometry was used. Fix \(S=\delta^{-\eta}\) with \(\eta=\epsilon^2\), and let \(\Delta:=S\delta\).

\noindent\textbf{Steps 1, 2, and 3.}
We construct refinements \((\Pcal_j,\M_j,\B_j)\), \(1\leq j\leq 4\), \((\Pcal_\Delta,\M_\Delta,\B_\Delta)\), and \(\B_{*}\) similarly to the first three steps in Theorem \ref{Thm_configuration}. In particular, let \(D_0\) be the fixed dilation from Lemma~\ref{lem:cinematic-duality}, and suppose \(|\mathbb P_4|\gtrapprox|\mathbb P_3|\) with
\begin{equation}\label{eq:uniform-dual-degree-cinematic}
  \#\{U\in\mathbb U_3: q\subseteq D_0\Gamma_U^\delta\}
  \sim W\qquad(q\in\mathbb Q_4),
\end{equation}
where \(\mathbb Q_4\) are the parameter cubes of \(\Pcal_4\).

\noindent\textbf{Step 4: High-low estimates.}
Similarly, consider the incidence estimates between the \(\delta\)-discretized balls \(\mathbb{B}_*\) and the \(\delta\times\delta^{1/2}\) rectangles \(\mathbb{U}_3\). Applying Theorem \ref{high_low_1} at scale \(\delta^{1/2}\), one of the following holds:
\begin{enumerate}[label=\textup{(\roman*)}]
  \item \label{high_1-cin} The high-frequency estimate gives
  \[
    \Inc^C(\M_3,\B_*)
      \lesssim_{\epsilon,C} \delta^{-1/4-\epsilon^2}|\M_3|^{1/2}|\B_*|^{1/2}.
  \]
  Equivalently, under the uniformity assumptions above,
   \begin{equation}\label{eq:planar-high-uniform-cin}
   N^2|\M_3|\lesssim_{\epsilon,C}\delta^{-1/2-2\epsilon^2}|\B_*|.
  \end{equation}
  \item \label{low_1-cin} The low-frequency estimate gives
  \begin{equation}\label{eq:planar-low-uniform-cin}
  I^C(\mathbb{U}_3,\mathbb{B}_*) \lesssim_{\epsilon,C} S^{-1+\epsilon^2} I_{\cap}^{C'}(\mathbb U_3,\mathbb B_*^{S\delta}).
  \end{equation}
\end{enumerate}

Next, apply Theorem~\ref{thm:cinematic-high-low} to \(\mathbb{Q}_4\) and \(\M_3\) using \eqref{eq:uniform-dual-degree-cinematic}. Put
\begin{equation}\label{eq:cin-IJ-using}
 \begin{split}
 \mathcal I&=\sum_{q\in\mathbb Q}
       \#\{U\in\mathbb U:q\subset D_0\Gamma_U^\delta\},\\
 \mathcal J&=\sum_{q\in\mathbb Q}
       \#\{U\in\mathbb U:
            q_{S\delta}\cap C'\Gamma_U^\delta
                     \ne\varnothing\}.
 \end{split}
\end{equation}
Then at least one of the following holds:
\begin{enumerate}[label=\textup{(\roman*)$'$}]
  \item\label{high_2-cin} The high-frequency estimate gives
  \begin{equation}\label{eq:dual-high-uniform-cin}
    W^4|\mathbb Q|
       \lesssim_{\epsilon,C,g}S^{O(1)}\delta^{-3/2}|\mathbb U|.
  \end{equation}
  \item\label{low_2-cin} The low-frequency estimate gives
  \begin{equation}\label{eq:dual-low-uniform-cin}
    \mathcal I
 \lesssim_{\epsilon,C,g}S^{-1+\epsilon^2} \mathcal  J.
  \end{equation}
\end{enumerate}

\medskip
\noindent\textbf{Both high-frequency estimates hold.}

Arguing as in Theorem \ref{Thm_configuration}, we obtain
\begin{align*}
      I^C_{\mathbb{U}_4}(\mathbb{P}_4,\mathbb{B}_*) \lesssim_{\epsilon,C,C_0,g} K^4 \delta^{1-\epsilon} |\mathbb{P}_4||\mathbb{B}_*|.
\end{align*}

\medskip
\noindent\textbf{The low-frequency estimate \ref{low_2-cin} holds.}

By the choice of \(N\) in Step 3 and
Lemma~\ref{lem:cinematic-duality},
\begin{align}
I^C_{\mathbb U_4}(\mathbb P_4,\mathbb B_*)
&\leq I^C_{\mathbb U_3}(\mathbb P_4,\mathbb B_*) \notag\\
&=
\sum_{p_q\in\mathbb P_4}
\sum_{\substack{U\in\mathbb U_3\\U\subset Cp_q}}
\#\{B\in\mathbb B_*:B\subset CU\} \notag\\
&\lesssim
N\sum_{q\in\mathbb Q_4}
\#\{U\in\mathbb U_3:q\subset D_0\Gamma_U^\delta\}
=N\mathcal I \notag\\
&\lesssim_{\epsilon,C,g}
S^{-1+\epsilon^2}N\mathcal J,
\label{eq:cinematic-low-first-transfer}
\end{align}
where the last inequality is
\eqref{eq:dual-low-uniform-cin}.

Similarly, consider
\begin{equation}\label{pair_2-cin}
   \Phi_g:=\#\left\{ (q,B)\in \mathbb Q_4\times \B_* : \exists U\in\mathbb{U}_3,\ \text{s.t.}\ B\subseteq CU,\ q_\Delta\cap C'\Gamma_U^\delta\neq \emptyset \right\}.
\end{equation}
On the one hand, for each \(U\in\mathbb{U}_3\) satisfying \eqref{pair_2-cin}, by the construction of \(\M'_\Delta\) in Step 2, there exists \(U_\Delta \in \M'_\Delta\) such that \(U \subseteq C_* U_\Delta\). After duality, \(\Gamma_U^\delta\subseteq C_* \Gamma^{\Delta}_{U_\Delta}\). Since \(q_\Delta\cap C' \Gamma_U^\delta\neq \emptyset\), Lemma \ref{lem:cinematic-jet-geometry} gives
\[
q_\Delta\subseteq\mathcal{N}_{10\Delta}(C'\Gamma_U^\delta)\subseteq C_gC'C_*\Gamma^{\Delta}_{U_\Delta}.
\]
By Lemma \ref{lem:cinematic-duality}, \(U_\Delta \subseteq D_1 p_q^\Delta\), where \(D_1\) is a fixed constant. Also, \(B_\Delta \subseteq D_1U_\Delta\). Hence
\begin{equation}\label{phi_l-cin}
    \Phi_g\lesssim I_{\M_{\Delta}}^{D_1}(\Pcal_4^{S\delta},\mathbb{B}^{S\delta}_*),
\end{equation}
where \(\Pcal_4^{S\delta}\) and \(\mathbb{B}^{S\delta}_*\) are counted with multiplicity.

On the other hand, if two boxes \(U_1, U_2\) satisfy condition \eqref{pair_2-cin} for the same pair \((q,B)\), then they must be \(D_3\)-comparable for some absolute constant \(D_3(C,C',g)\). Indeed, if \(U_1,U_2\) satisfy the above, then \(B\subseteq CU_1\cap CU_2\) implies
\[
|t_{U_j}-t_B|\lesssim C\sqrt\delta,\quad
 |y_B-y_{U_j}-s_{U_j}(t_B-t_{U_j})|\lesssim C\delta,\quad j=1,2,
\]
where \((t_B,y_B)\) is the center of \(B\). By Lemma \ref{lem:cinematic-jet-geometry},
\[
q_\Delta \subseteq \mathcal{N}_{2\Delta} (C'\Gamma_{U_j}^\delta) \subseteq C''\Gamma_{U_j}^{\Delta,\delta^{1/2}},\quad j=1,2.
\]
Hence, there exists \(z\in q_\Delta\) such that
\[
|g_t(z,t_{U_j})-s_{U_j}|\lesssim C''\delta^{1/2}, \quad j=1,2.
\]
Thus
\[
|t_{U_1}-t_{U_2}|,\ |s_{U_1}-s_{U_2}|\lesssim_{g,C,C''} \delta^{1/2}, \quad |y_{U_{1}}-y_{U_2}-s_{U_2} (t_{U_1}-t_{U_2})|\lesssim_{g,C,C''} \delta,
\]
so \(U_1,U_2\) are \(D_3\)-comparable. Therefore,
\begin{equation}\label{Phi_g-cin}
    \Phi_g \sim_{C} N
    \mathcal{J},
\end{equation}
and
\begin{equation}\label{eq:dual-low-contraction-cin}
  I^C_{\mathbb U_4}(\mathbb P_4,\mathbb B_*)
  \lesssim_{\epsilon,C_0,C}
   S^{-1+\epsilon^2}I_{\M_{\Delta}}^{D_1}(\Pcal_4^{S\delta},\mathbb{B}^{S\delta}_*).
\end{equation}
The conclusion then follows by induction on \(S\delta\), as in Theorem \ref{Thm_configuration}.

\noindent\textbf{The low-frequency estimate \ref{low_1-cin} holds.}

Similarly,
\begin{align}
    I^C_{\mathbb{U}_4}(\mathbb{P}_4,\mathbb{B}_*)  \lesssim_{\epsilon,C} S^{-1+\epsilon^2} MI_{\cap}^{C'}(\mathbb U_3,\mathbb B_*^{S\delta}).
\end{align}
Consider
\begin{equation}\label{pair_1-cin}
   \Psi_g:=\#\{ (p,B)\in \Pcal_2\times \B_* : \exists U\in\mathbb{U}_3,\ U\subseteq Cp,\ B_\Delta\cap C' U\neq \emptyset \}.
\end{equation}
On the one hand, for each \(U\in\mathbb{U}'_3\) satisfying \eqref{pair_1-cin}, choose \(U_\Delta \in \M'_\Delta\) such that \(U \subseteq C_* U_\Delta\). Since \(B_\Delta\cap C'U\neq \emptyset\) and \(U\subseteq Cp\), the Taylor formula implies that there exists a constant \(D_2\) (depending on \(C,C',C_*,g\)) such that
\[
B_\Delta\subseteq D_2  U_\Delta,\qquad U_\Delta \subseteq D_2 p^\Delta.
\]
Hence,
\begin{equation}
    \label{Psi_l-cin}
    \Psi_g \lesssim I_{\M_{\Delta}}^{D_2}(\Pcal_2^{S\delta},\mathbb{B}^{S\delta}_*).
\end{equation}

On the other hand, similarly, if two boxes \(U_1,U_2\) satisfy \eqref{pair_1-cin} for the same \((p,B)\), then they must be \(D_4\)-comparable for some \(D_4\). Suppose that \(U_1, U_2 \subseteq Cp\) are distinct. Then they are contained in \(\delta^{1/2}\)-charts \(\omega_1, \omega_2\) on the \(\delta\)-cinematic curve \(p\), which are separated by at least \(C\delta^{1/2}\). It follows that
\[
B_\Delta \subseteq \mathcal{N}_{C'\Delta}(\omega_1) \cap \mathcal{N}_{C'\Delta}(\omega_2),
\]
which yields a contradiction. Thus
\begin{equation}
    \label{Psi_g-cin}
   \Psi_g \sim_C MI_{\cap}^{C'}(\mathbb U_3,\mathbb B_*^{S\delta}).
\end{equation}
Combining the preceding estimates,
\begin{equation}
  I^C_{\mathbb U_4}(\mathbb P_4,\mathbb B_*)
  \lesssim_{\epsilon,C_0,C}
   S^{-1+\epsilon^2}I_{\M_{\Delta}}^{D_2}(\Pcal_2^{S\delta},\mathbb{B}^{S\delta}_*).
\end{equation}
Thus the induction closes similarly.

\end{proof}

\section {from discretized estimates to Hausdorff dimension bounds}
\label{section 7}
The following proposition is a straightforward adaptation of Proposition 3.1 in \cite{GreenHarrisOuRenTammen2025Circular}, which converts incidence bounds into lower bounds for the Hausdorff dimension of Furstenberg sets.

\begin{proposition}
     Theorem \ref{Thm_discretized} implies Theorem \ref{main_thm_Hausdorff_form} and Theorem \ref{thm:cinematic-hausdorff}.
\end{proposition}
\begin{proof}
Let \(F\) be a circular \((u,v)\)-Furstenberg set with
\(2u+v\geq3\), and write \(d=\dim_{\mathrm H}F\). Fix \(\tau>0\) and let
\(\alpha>d\). Apply Theorem~\ref{Thm_discretized} with
\(\epsilon=\tau\), and let \(\eta>0\) and \(\delta_0>0\) be the
corresponding constants. Choose
\[
0<\rho<\frac{1}{10}\min\{u,v,\eta\},
\qquad u_0:=u-\rho,\qquad v_0:=v-\rho.
\]

Let \(E\) be a family of circles witnessing the Furstenberg property. By scaling and countable sub-additivity of the Hausdorff content, we can assume the circles defining the Furstenberg property for F have radii between 1 and 2.

Since \(u_0<u\) and \(v_0<v\), a standard countable localization in the
circle parameters and the fiber contents(also see the proof of  Lemma 3.3 in \cite{hera2022improved} ) , there exists a constant \(c>0\) and $E_0 \subseteq E$ such that
\[
\mathcal H^{v_0}_\infty(E_0)\geq c,
\qquad
\mathcal H^{u_0}_\infty(F\cap C_z)\geq c
\quad (z\in E_0).
\]
Here \(C_z\) is the circle with parameter \(z\). All circles in \(E_0\)
lie in a fixed bounded planar set \(\Omega\). Replacing \(F\) by
\(F\cap\Omega\) does not change any of the intersections above, so we
may assume that \(F\) is bounded.

Let \(k_0\) be arbitrarily large, with \(2^{-k_0}<\delta_0\). Since
\(\alpha>d\), the bounded part of \(F\) under consideration has a
dyadic cover
\[
F\subseteq
\bigcup_{k\geq k_0}\bigcup_{D\in\mathcal T_k}D,
\qquad
\sum_{k\geq k_0}2^{-k\alpha}|\mathcal T_k|\leq 1,
\]
where the elements of \(\mathcal T_k\) are dyadic squares of side length
\(2^{-k}\). Similarly, by the subadditivity of Hausdorff content, for each
\(z\in E_0\) there is an index \(k\geq k_0\) such that
\[
\mathcal H^{u_0}_\infty
\left(C_z\cap\bigcup_{D\in\mathcal T_k}D\right)
\gtrsim k^{-2}.
\]
Assign one such index $k$ to each \(z\), and denote the resulting level set
by \(E_k\). Thus \(E_0=\bigcup_{k\geq k_0}E_k\). A second application of
subadditivity gives an index \(k\geq k_0\) for which
\[
\mathcal H^{v_0}_\infty(E_k)\gtrsim k^{-2}.
\]
Set \(\delta=2^{-k}\). We use the following standard dyadic consequence of Frostman's lemma.
If \(A\subset\mathbb R^n\) is bounded,
\(\mathcal H^s_\infty(A)\geq a>0\), \(\delta=2^{-k}\), and
\(\mathcal Q\) is a family of dyadic \(\delta\)-cubes covering \(A\),
then \(\mathcal Q\) contains a subfamily
\[
\mathcal Q'\subseteq
\{Q\in\mathcal Q:Q\cap A\neq\emptyset\}
\]
which is a
\((\delta,s,a^{-1}k^{O(1)})\)-set. Indeed, we can get \(\mathcal Q'\) by applying Frostman's lemma to
\(A\) and pigeonholing the masses of the cubes in \(\mathcal Q\). (see Lemma 2.13\cite{F_ssler_2014})

Apply the dyadic Frostman extraction above to the family of dyadic cubes of side length \(\delta\) meeting \(E_k\). The resulting family
\(\mathbb P\) of discretized circles is a
\((\delta,v_0,k^{O(1)})\)-set. For each \(p\in\mathbb P\), choose \(z_p\in \hat{p}\) . Applying
the same extraction to
\[
C_{z_p}\cap\bigcup_{D\in\mathcal T_k}D
\]
gives a subfamily \(\mathbb B(p)\subseteq\mathcal T_k\) which is a
\((\delta,u_0,k^{O(1)})\)-set. If \(B\in\mathbb B(p)\), then
\(B\cap C_{z_p}\neq\emptyset\), \(\operatorname{diam}(B)\lesssim\delta\),
and \(z_p\in\widehat p\). Hence, for an absolute constant \(C_0\),
\[
B\subseteq\mathcal N_{O(\delta)}(C_{z_p})\subseteq C_0p.
\]
Consequently, with
\(
\mathbb B:=\bigcup_{p\in\mathbb P}\mathbb B(p),
\)
the pair \((\mathbb B,\mathbb P)_\delta\) is a circular Furstenberg
shading model.

For every \(\delta\leq r\leq1\),
\[
r^{u_0}\leq\delta^{-\rho}r^u,
\qquad
r^{v_0}\leq\delta^{-\rho}r^v.
\]
It follows that \(\mathbb P\) is a
\((\delta,v,k^{O(1)}\delta^{-\rho})\)-set and every
\(\mathbb B(p)\) is a
\((\delta,u,k^{O(1)}\delta^{-\rho})\)-set. Since \(\rho<\eta\), for all
sufficiently large \(k_0\), and hence for every \(k\geq k_0\),
\[
k^{O(1)}\delta^{-\rho}\leq\delta^{-\eta}.
\]
Thus \((\mathbb B,\mathbb P)_\delta\) satisfies the hypotheses of
Theorem~\ref{Thm_discretized}. Since \(2u+v\geq3\), we obtain
\[
|\mathbb B|\gtrsim_{\tau,C_0}
\delta^\tau\delta^{-u-1}.
\]
On the other hand, \(\mathbb B\subseteq\mathcal T_k\), and the choice of
the cover gives
\[
|\mathbb B|\leq|\mathcal T_k|\leq\delta^{-\alpha}.
\]
Since \(k_0\) was arbitrary, \(\delta\) can be made arbitrarily small.
Therefore \(\alpha\geq u+1-\tau\). Letting \(\alpha\downarrow d\) and
then \(\tau\downarrow0\) yields
\[
\dim_{\mathrm H}F\geq u+1,
\]
as required.

The cinematic case is similar. After suitable localization, we can obtain a constant \(c>0\) and $E_0 \subseteq E$ such that
\[
\mathcal H^{v_0}_\infty(E_0)\geq c,
\qquad
\mathcal H^{u_0}_\infty(F\cap \Sigma_z)\geq c
\quad (z\in E_0).
\]
Then the problem can be reduced to discretized estimate by applications of Frostman's lemma.
\end{proof}

\bibliographystyle{amsalpha}
\bibliography{references}

\end{document}